\documentclass{compositio}
\usepackage{amsmath}
\usepackage{amsthm}
\usepackage{amssymb}
\usepackage{amsfonts}
\usepackage{amsxtra}
\usepackage{epsfig}
\usepackage{verbatim}
\usepackage{enumerate}
\usepackage{amscd}
\theoremstyle{plain}
\newtheorem{thm}{Theorem}
\newtheorem{prop}[thm]{Proposition}
\newtheorem{lem}[thm]{Lemma}
\newtheorem{cor}[thm]{Corollary}

\theoremstyle{definition}

\newtheorem{conjecture}[thm]{Conjecture}

\theoremstyle{remark}
\newtheorem{example}[thm]{Example}
\newtheorem*{rem}{Remark}

\include{header}
\newcommand{\Real}{\mathbf R}
\newcommand{\C}{\mathbf C}
\newcommand{\N}{\mathbf N}
\newcommand{\Q}{\mathbf Q}
\newcommand{\Z}{\mathbf Z}
\newcommand{\A} {\mathbf{A}}
\newcommand{\h} {\mathbf{H}_3}
\newcommand{\oh} {\overline{\mathbf{H}}_3}

\newcommand{\Ox}{\mathcal{R}}
\newcommand{\Fx}{L}
\newcommand{\Tm}{\mathbf{T}_{\mathfrak{m}}}
\newcommand{\Eis}{{\rm Eis}(\phi)}
\begin{document}

\title[Eisenstein ideal for imaginary quadratic fields and Bloch-Kato]{An Eisenstein ideal for imaginary quadratic fields and the
Bloch-Kato conjecture for Hecke characters}

\author{Tobias Berger}
\email{t.berger@dpmms.cam.ac.uk}
\address{Department of Pure Mathematics and Mathematical Statistics, Centre for Mathematical Sciences, University of Cambridge,
Cambridge CB3 0WB, United Kingdom}



\begin{abstract}
For certain algebraic Hecke characters $\chi$ of an imaginary
quadratic field $F$ we define an Eisenstein ideal in a $p$-adic
Hecke algebra acting on cuspidal automorphic forms of ${\rm
GL}_2/F$. By finding congruences between Eisenstein cohomology
classes (in the sense of G. Harder) and cuspidal classes we prove a
lower bound for the index of the Eisenstein ideal in the Hecke
algebra in terms of the special $L$-value $L(0,\chi)$. We further
prove that its index is bounded from above by the order of the
Selmer group of the $p$-adic Galois character associated to
$\chi^{-1}$. This uses the work of R. Taylor \textit{et al.} on
attaching Galois representations to cuspforms of ${\rm GL}_2/F$.
Together these results imply a lower bound for the size of the
Selmer group in terms of $L(0,\chi)$, coinciding with the value
given by the Bloch-Kato conjecture.
\end{abstract}

\maketitle

\section*{Introduction}
The aim of this work is to demonstrate the use of Eisenstein
cohomology, as developed by Harder, in constructing elements of
Selmer groups for Hecke characters of an imaginary quadratic field
$F$. The strategy of first finding congruences between Eisenstein
series and cuspforms and then using the Galois representations
attached to the cuspforms to prove lower bounds on the size of
Selmer groups goes back to Ribet \cite{Ri}, and has been applied and
generalized in \cite{Wi} in his proof of the Iwasawa main conjecture
for characters over totally real fields, and more recently by
\cite{SU} and \cite{BeCh} amongst others. These all used for the
congruences integral structures coming from algebraic geometry. In
our case the symmetric space associated to ${\rm GL}_2/F$ is not
hermitian and we therefore use the integral structure arising from
Betti cohomology. This alternative, more general approach was
outlined for ${\rm GL}_{2}/\Q$ in \cite{HaPi}.

In \cite{HaGL2} G. Harder constructed Eisenstein cohomology as a
complement to the cuspidal cohomology for the groups ${\rm GL}_2$
over number fields and proved that this decomposition respects the
rational structure of group cohomology. The case interesting for
congruences is when this decomposition is not integral, i.e., when
there exists an Eisenstein class with integral restriction to the
boundary that has a denominator. In \cite{Be06} such classes were
constructed for imaginary quadratic fields and their denominator was
bounded from below by the special $L$-value of a Hecke character. In
\S \ref{s5} we give a general set-up in which cohomological
congruences between Eisenstein and cuspidal classes can be proven
(Proposition \ref{Eisideal}) and then apply this to the classes
constructed in \cite{Be06}. The result can be expressed as a lower
bound for the index of the Eisenstein ideal of the title in a Hecke
algebra in terms of the special $L$-value. The main obstacle to
obtaining a congruence from the results in \cite{Be06} is the
occurrence of torsion in higher degree cohomology, which is not very
well understood (see the discussion in \S \ref{tproblem}). However,
we manage to solve this ``torsion problem" for unramified characters
in \S \ref{s6}.

The other main ingredient are the Galois representations associated
to cohomological cuspidal automorphic forms, constructed by R.
Taylor \textit{et al.} by means of a lifting to the symplectic group
in \cite{T2}. Assuming the existence of congruences between an
Eisenstein series and cuspforms we use these representations in \S
\ref{SelEis} to construct elements in the Selmer group of a Galois
character. In fact, we prove that its size is bounded from below by
the index of the Eisenstein ideal.

These two results are combined in \S \ref{s7} to prove a lower
bound on the size of Selmer groups of Hecke characters of an
imaginary quadratic field in terms of a special $L$-value,
coinciding with the value given by the Bloch-Kato conjecture.

To give a more precise account, let $p>3$ be a prime unramified in
the extension $F/\Q$ and let $\mathfrak{p}$ be a prime of $F$
dividing $(p)$. Fix embeddings $\overline F \hookrightarrow
\overline F_{\mathfrak{p}} \hookrightarrow \C$. Let $\phi_1,
\phi_2:F^* \backslash \A_F^* \to \C^*$ be two Hecke characters of
infinity type $z$ and $z^{-1}$, respectively. Let $\Ox$ be the ring
of integers in a sufficiently large finite extension of
$F_{\mathfrak{p}}$. Let $\mathbf{T}$ be the $\Ox$-algebra generated
by Hecke operators acting on cuspidal automorphic forms of ${\rm
GL}_2/F$. For $\phi=(\phi_1, \phi_2)$ we define in \S \ref{SelEis}
an Eisenstein ideal $\mathbf{I}_{\phi}$ in $\mathbf{T}$. Following
previous work of Wiles and Urban we construct elements in the Selmer
group of $\chi_{\mathfrak{p}} \epsilon$, where $\chi_{\mathfrak{p}}$
is the $p$-adic Galois characters associated to
$\chi:=\phi_1/\phi_2$ and $\epsilon$ is the $p$-adic cyclotomic
character. We obtain a lower bound on the size of the Selmer group
in terms of that of the congruence module
$\mathbf{T}/\mathbf{I}_{\phi}$. A complication that arises in the
application of Taylor's theorem is that we need to work with
cuspforms with cyclotomic central character. This is achieved by a
twisting argument (see Lemma \ref{anticycchar}).

To prove the lower bound on the congruence module in terms of the
special $L$-value (the first step described above), we use the
Eisenstein cohomology class ${\rm Eis}(\phi)$ constructed in
\cite{Be06} in the cohomology of a symmetric space $S$ associated to
${\rm GL}_2/F$. The class is an eigenvector for the Hecke operators
at almost all places with eigenvalues corresponding to the
generators of $\mathbf{I}_{\phi}$, and its restriction to the
boundary of the Borel-Serre compactification of $S$ is integral. The
main result of \cite{Be06}, which we recall in \S 5, is that the
denominator $\delta$ of ${\rm Eis}(\phi) \in H^1(S,\overline
F_{\mathfrak{p}})$ is bounded from below by $L^{\rm alg}(0,\chi)$.
As mentioned above, Proposition \ref{Eisideal} gives a general setup
for cohomological congruences. It implies the existence of a
cuspidal cohomology class congruent to $\delta \cdot {\rm
Eis}(\phi)$ modulo the $L$-value supposing that there exists an
integral cohomology class with the same restriction to the boundary
as ${\rm Eis}(\phi)$. The latter can be replaced by the assumption
that $H^2_c(S, \Ox)_{\rm torsion}=0$, and this result is given in
Theorem \ref{condthm}. In \S \ref{s6} we prove that the original
hypothesis is satisfied for unramified $\chi$, avoiding the issue of
torsion freeness. We achieve this by a careful analysis of the
restriction map to the boundary $\partial \overline S$ of the
Borel-Serre compactification. Starting with a group cohomological
result for $\mathrm{SL}_2(\mathcal{O})$ due to Serre \cite{Se}
(which we extend to all maximal arithmetic subgroups of
$\mathrm{SL}_2(F)$) we define an involution on $H^1(\partial
\overline S, \Ox)$ such that the restriction map
$$H^1(S, \Ox) \overset{\mathrm{res}}{\twoheadrightarrow}
H^1(\partial \overline S, \Ox)^-,$$ surjects onto the
$-1$-eigenspace.  We apply the resulting criterion to
$\mathrm{res}(\mathrm{Eis}(\phi))$ to deduce the existence of a
lift to $H^1(S, \Ox)$.

Note that the restriction to constant coefficient systems and
therefore weight 2 automorphic forms is important only for \S
\ref{s6}. It was applied throughout to simplify the exposition. In
particular, the results of Theorems \ref{thmdenom} and \ref{condthm}
extend (for split $p$) to characters $\chi$ with infinity type
$z^{m+2} \overline z^{-m}$ for $m \in \N_{\geq 0}$. See \cite{Be06}
for the necessary modifications and results.

Combining the two steps we obtain in \S \ref{s7} a lower bound for
the size of the Selmer group of $\chi_{\mathfrak{p}} \epsilon$ in
terms of $L^{\rm alg}(0,\chi)$ and  relate this result to the
Bloch-Kato conjecture. This conjecture has been proven in our case
(at least for class number $1$) starting from the Main Conjecture of
Iwasawa theory for imaginary quadratic fields (see \cite{Han},
\cite{Guo93}).  Similar results have also been obtained by Hida in
\cite{Hi82} for split primes $p$ and $\overline \chi=\chi^c$ using
congruences of classical elliptic modular forms between CM and
non-CM forms. We seem to recover base changes of his congruences in
this case (see \S \ref{discussion}).

However, our method of constructing elements in Selmer groups using
cohomological congruences is very different and should be more
widely applicable. The analytic theory of Eisenstein cohomology has
been developed for many groups, and rationality results are known,
e.g. for ${\rm GL}_n$ by the work of \cite{FS}. Our hope is that the
method presented here generalizes to these higher rank groups.

To conclude, we want to mention two related results. In \cite{F}
congruences involving degree 2 Eisenstein cohomology classes for
imaginary quadratic fields were constructed but only the $L$-value
of the quadratic character associated to $F/\Q$ was considered. The
torsion problem we encounter does not occur for degree 2, but the
treatment of the denominator of the Eisenstein classes is more
difficult. For cases of the Bloch-Kato conjecture when the Selmer
groups are infinite see \cite{BeCh}. Their method is similar to ours
in that they use congruences between Eisenstein series and
cuspforms, however, they work with $p$-adic families on ${\rm U}(3)$
and do not use Eisenstein cohomology.

\section{Notation and Definitions}
\subsection{General notation} Let $F/\Q$ be an
imaginary quadratic extension and $d_F$ its absolute discriminant.
Denote the classgroup by ${\rm Cl}(F)$ and the ray class group
modulo a fractional ideal $\mathfrak{m}$ by ${\rm
Cl}_{\mathfrak{m}}(F)$. For a place $v$ of $F$ let $F_v$ be the
completion of $F$ at $v$. We write $\mathcal{O}$ for the ring of
integers of $F$, $\mathcal{O}_v$ for the closure of $\mathcal{O}$ in
$F_v$, $\mathfrak{P}_v$ for the maximal ideal of $\mathcal{O}_v$,
$\pi_v$ for a uniformizer of $F_v$, and $\hat{\mathcal{O}}$ for
$\prod_{v \, {\rm finite}} \mathcal{O}_v$. We use the notations $\A,
\A_f$ and $\A_F, \A_{F,f}$ for the adeles and finite adeles of $\Q$
and $F$, respectively, and write $\A^*$ and $\A_F^*$ for the
respective group of ideles. Let $p$ be a prime of $\Z$ that does not
ramify in $F$, and let $\mathfrak{p} \subset \mathcal{O}$ be a prime
dividing $(p)$.

Denote by $G_F$ the absolute Galois group of $F$. For $\Sigma$ a
finite set of places of $F$ let $G_{\Sigma}$ be the Galois group of
the maximal extension of $F$ unramified at all places not in
$\Sigma$. We fix an embedding $\overline F \hookrightarrow \overline
F_v$ for each place $v$ of $F$. Denote the corresponding
decomposition and inertia groups by $G_v$ and $I_v$, respectively.
Let $g_v=G_v/I_v$ be the Galois group of the maximal unramified
extension of $F_v$. For each finite place $v$ we also fix an
embedding $\overline F_v \hookrightarrow \C$ that is compatible with
the fixed embeddings $i_v: \overline F \hookrightarrow \overline
F_v$ and $i_{\infty}: \overline F \hookrightarrow \C(=\overline
F_{\infty})$. For a topological $G_F$-module (resp. $G_v$-module)
$M$ write $H^1(F,M)$ for the continuous Galois cohomology group
$H^1(G_F,M)$, and $H^1(F_v,M)$ for $H^1(G_v,M)$.

\subsection{Hecke characters} \label{s2.2} A Hecke character of $F$ is a continuous group homomorphism
 $\lambda: F^* \backslash \A_F^* \to \C^*$. Such a character corresponds
 uniquely to a character on ideals prime to the conductor, which we also denote by $\lambda$.
 Define the character $\lambda^c$ by
$\lambda^c(x)=\lambda(\overline x)$.
\begin{lem}[(Lemma 3.1 of \cite{Be})] \label{unanticyc}
If $\lambda$ is an unramified Hecke character then
$\lambda^c=\overline \lambda$. \hspace{\fill} \qedsymbol
\end{lem}
 For Hecke characters $\lambda$ of type $(A_0)$, i.e., with
infinity type $\lambda_{\infty}(z)=z^m \overline z^n$ with $m, n
\in \Z$ we define (following Weil) a $p$-adic Galois character
$$\lambda_{\mathfrak{p}}:G_F \to \overline F_{\mathfrak{p}}^*$$ associated to $\lambda$ by the following
rule: For a finite place $v$ not dividing $p$ or the conductor of
$\lambda$, put $\lambda_{\mathfrak{p}}({\rm
Frob}_v)=i_{\mathfrak{p}}(i_{\infty}^{-1}(\lambda(\pi_v)))$ where
${\rm Frob}_v$ is the \textit{arithmetic} Frobenius at $v$. It
takes values in the integer ring of a finite extension of
$F_{\mathfrak{p}}$.

Let $\epsilon: G_F \to \Z_p^*$ be the $p$-adic cyclotomic
character defined by the action of $G_F$ on the $p$-power roots of
unity: $g.\xi=\xi^{\epsilon(g)}$ for $\xi$ with $\xi^{p^m}=1$ for
some $m$. Our convention is that the Hodge-Tate weight of
$\epsilon$ at $\mathfrak{p}$ is $1$.

Write $L(0,\lambda)$ for the Hecke $L$-function of $\lambda$. Let
$\lambda$ a Hecke character of infinity type $z^a \left(
\frac{z}{\overline z} \right )^{b}$ with conductor prime to $p$.
Assume $a,b \in \Z$ and $a>0$ and ${b} \geq 0$. Put
$$L^{\mathrm{alg}}(0,\lambda):= \Omega^{-a-2{b}} \left(
\frac{2\pi}{\sqrt{d_F}}\right)^{b} \Gamma(a+{b})\cdot
L(0,\lambda).$$ In most cases, this normalization is integral,
i.e., lies in the integer ring of a finite extension of
$F_{\mathfrak{p}}$. See \cite{Be06} Theorem 3 for the exact
statement. Put $$L^{\rm int}(0, \lambda)=
  \begin{cases}
    L^{\mathrm{alg}}(0,\lambda) & \text{if } {\rm val}_p (L^{\mathrm{alg}}(0,\lambda)) \geq 0 \\
    1 & \text{otherwise}.
  \end{cases}
$$

\subsection{Selmer groups}
Let $\rho:G_F \to \mathcal{R}^*$ be a continuous Galois character
taking values in the ring of integers $\mathcal{R}$ of a finite
extension $L$ of $F_{\mathfrak{p}}$. Write
$\mathfrak{m}_{\mathcal{R}}$ for its maximal ideal and put
$\mathcal{R}^{\vee}=L/\mathcal{R}$. Let $\mathcal{R}_{\rho}$,
$L_{\rho}$, and
$W_{\rho}=L_{\rho}/\mathcal{R}_{\rho}=\mathcal{R}_{\rho}
\otimes_{\mathcal{R}} \mathcal{R}^{\vee}$ be the free rank one
modules on which $G_F$ acts via $\rho$.

Following Bloch and Kato \textit{et al.} we define the following
Selmer groups: Let $$H^1_f(F_v, L_{\rho})=
  \begin{cases}
    {\rm ker}(H^1(F_v, L_{\rho}) \to H^1(I_v, L_{\rho})) & \text{ for } v \nmid p, \\
    {\rm ker}(H^1(F_v, L_{\rho}) \to H^1(F_v, B_{\rm cris} \otimes L_{\rho})) & \text{ for } v \mid
    p,
  \end{cases}$$ where $B_{\rm cris}$ denotes Fontaine's ring of  $p$-adic
  periods. Put $$H^1_f(F_v, W_{\rho})= {\rm
im}(H^1_f(F_v, L_{\rho}) \to H^1(F_v, W_{\rho})).$$ For a finite set
of places $\Sigma$ of $F$ define
$${\rm Sel}^{\Sigma}(F,\rho)={\rm ker}\left( H^1(F,W_{\rho}) \to \prod_{v
\notin \Sigma} \frac{H^1(F_v, W_{\rho})}{H^1_f(F_v,W_{\rho})}
\right).$$ We write ${\rm Sel}(F,\rho)$ for ${\rm
Sel}^{\emptyset}(F,\rho)$.

If $p$ splits in $F/\Q$ and $\rho=\lambda_{\mathfrak{p}}$ for a
Hecke character $\lambda$ of infinity type $z^a \overline z^b$ with
$a,b \in \Z$ (``ordinary case") then we define
$$F_{\mathfrak{p}}^+ L_{\rho}=
  \begin{cases}
    L_{\rho} & \text{ if }a<0 \text{ (i.e., HT-wt of } \rho >0), \\
    \{0\} & \text{ if } a \geq 0 \text{ (i.e., HT-wt of } \rho \leq 0)
  \end{cases}$$
and $$F_{\overline{\mathfrak{p}}}^+ L_{\rho}=
  \begin{cases}
    L_{\rho} & \text{ if }b<0 , \\
    \{0\} & \text{ if } b \geq 0.
  \end{cases}$$
In the ordinary case we have $H^1_f(F_v, L_{\rho})=H^1(F_v, F^+_v
L_{\rho})$ for $v \mid p$ (see \cite{Guo96} p.361, \cite{Fl} Lemma
2).

\begin{lem} \label{tameinertia}
Let $\rho$ be unramified at $v \nmid p$. If $\rho({\rm Frob}_v)
\not \equiv \epsilon({\rm Frob}_v) \mod{p}$ then
$${\rm Sel}^{\Sigma}(F,\rho) =
{\rm Sel}^{\Sigma \backslash \{v\}}(F,\rho).$$
\end{lem}

\begin{proof}
By definition ${\rm Sel}^{\Sigma \backslash \{v\}}(F,\rho) \subset
{\rm Sel}^{\Sigma}(F,\rho)$ for any $v$. For places $v$ as in the
lemma we have
$$H^1_f(F_v, W_{\rho})=
{\rm ker}(H^1(F_v, W_{\rho}) \to H^1(I_v, W_{\rho})^{g_v}).$$ It is
clear that $H^1(I_v, W_{\rho})^{g_v}={\rm Hom}_{g_v}(I_v^{\rm
tame},W_{\rho})={\rm Hom}_{g_v}(I_v^{\rm
tame},W_{\rho}[\mathfrak{m}_{\mathcal{R}}^n])$ for some $n$. By our
assumption therefore $H^1(I_v, W_{\rho})^{g_v}=0$ since ${\rm
Frob}_v$ acts on $I_v^{\rm tame}$ by $\epsilon({\rm Frob}_v)$.
\end{proof}

\subsection{Cuspidal automorphic representations}
We refer to \cite{U95} \S 3.1 as a reference for the following: For
$K_f=\prod_v K_v \subset G(\A_f)$ a compact open subgroup, denote by
$S_2(K_f)$ the space of cuspidal automorphic forms of
$\mathrm{GL}_2(F)$ of weight 2, right-invariant under $K_f$. For
$\omega$ a finite order Hecke character write $S_2(K_f, \omega)$ for
the forms  with central character $\omega$. This is isomorphic as a
$G(\A_f)$-module to $\bigoplus \pi_f^{K_f}$ for automorphic
representations $\pi$ of a certain infinity type (see Theorem
\ref{taylor} below) with central character $\omega$. For $g\in
G(\A_f)$ we have the Hecke action of $[K_f g K_f]$ on $S_2(K_f)$ and
$S_2(K_f, \omega)$. For places $v$ with $K_v={\rm
GL}_2(\mathcal{O}_v)$ we define $T_v=[K_f
\begin{pmatrix} \pi_v &0\\0&1 \end{pmatrix} K_f]$.

\subsection{Cohomology of symmetric space} \label{s2.5}
Let $G= {\rm Res}_{F/\Q} {\rm GL}_2$ and $B$ the restriction of
scalars of the Borel subgroup of upper triangular matrices. For any
$\Q$-algebra $R$ we consider a pair of characters $\phi=(\phi_1,
\phi_2)$ of $R^* \times R^*$ as character of $B(R)$ by defining
$\phi(\begin{pmatrix} a & b\\ 0 & d \end{pmatrix})=\phi_1(a)
\phi_2(b)$. Put $K_{\infty}=U(2) \cdot \C^* \subset G(\Real)$. For
an open compact subgroup $K_f \subset G(\A_f)$ we define the adelic
symmetric space
$$S_{K_f}=G(\Q)\backslash G(\A) /K_{\infty} K_f.$$
Note that $S_{K_f}$ has several connected components. In fact,
strong approximation implies that the fibers of the determinant
map
$$S_{K_f} \twoheadrightarrow \pi_0(K_f):=\A_{F,f}^*/{\rm det}(K_f)F^*$$ are
connected. Any $\gamma \in G(\A_f)$ gives rise to an injection
$$G_{\infty} \hookrightarrow G(\A)$$ $$g_{\infty} \mapsto
(g_{\infty},\gamma)$$ and, after taking quotients, to a component
$\Gamma_{\gamma}\backslash G_{\infty}/K_{\infty} \to S_{K_f},$
where $$\Gamma_{\gamma}:= G(\Q) \cap \gamma K_f \gamma^{-1}.$$
This component is the fiber over ${\rm det}(\gamma)$. Choosing a
system of representatives for $\pi_0(K_f)$ we therefore have
$$ S_{K_f} \cong \coprod_{[{\rm det}(\gamma)] \in \pi_0(K_f)}
\Gamma_{\gamma} \backslash \h,$$ where $G_{\infty}/K_{\infty}$ has
been identified with three-dimensional hyperbolic space
$\h=\Real_{>0} \times \C$.

We denote the Borel-Serre compactification of $S_{K_f}$ by
$\overline S_{K_f}$ and write $\partial \overline S_{K_f}$ for its
boundary. The Borel-Serre compactification $\overline S_{K_f}$ is
given by the union of the compactifications of its connected
components. For any arithmetic subgroup $\Gamma \subset G(\Q)$,
the boundary of the Borel-Serre compactification of $\Gamma
\backslash \h$, denoted by $\partial(\Gamma \backslash \oh)$, is
homotopy equivalent to
\begin{equation} \label{bdry} \coprod_{[\eta] \in
\mathbf{P}^1(F)/\Gamma} \Gamma_{B^{\eta}}\backslash \h,
\end{equation} where we identify $\mathbf{P}^1(F)=B(\Q) \backslash G(\Q)$, take $\eta \in G(\Q)$,
 and put $\Gamma_{B^{\eta}}=\Gamma \cap \eta^{-1}B(\Q)\eta$.

For $X \subset \overline S_{K_f}$ and $R$ an $\mathcal{O}$-algebra
we denote by $H^i(X,R)$ (resp. $H^i_c(X,R)$)  the $i$-th (Betti)
cohomology group  (resp. with compact support), and the interior
cohomology, i.e., the image of $H^i_c(X,R)$ in $H^i(X,R)$, by
$H^i_!(X,R)$.

There is a Hecke action of double cosets $[K_f g K_f]$ for $g\in
G(\A_f)$ on these cohomology groups (see \cite{U98} \S 1.4.4 for
the definition). We put $T_{\pi_v}=[K_f \begin{pmatrix} \pi_v &
0\\ 0 & 1 \end{pmatrix}K_f]$ and $S_{\pi_v}=[K_f
\begin{pmatrix} \pi_v & 0\\ 0 & \pi_v
\end{pmatrix}K_f]$.

The connection between cohomology and cuspidal automorphic forms
is given by the Eichler-Shimura-Harder isomorphism (in this
special case see \cite{U98} Theorem 1.5.1): For any compact open
subgroup $K_f \subset G(\A_f)$ we have
\begin{equation} \label{ESH}
S_2(K_f) \overset{\sim}{\to} H^1_!(S_{K_f}, \C)
\end{equation} and the isomorphism is Hecke-equivariant.

One knows (see for example, \cite{Be06} Proposition 4) that for
any $\mathcal{O}[\frac{1}{6}]$-algebra $R$ there is a natural
$R$-functorial isomorphism \begin{equation} \label{shgp}
H^1(\Gamma \backslash \h, R) \cong H^1(\Gamma, R),
\end{equation} where the group cohomology $H^1(\Gamma, R)$ is just given by
${\rm Hom}(\Gamma,R)$.

\subsection{Galois representations associated to cuspforms for
imaginary quadratic fields} \label{Galreps} Combining the work of
Taylor, Harris, and Soudry with results of Friedberg-Hoffstein and
Laumon/Weissauer, one can show the following (see \cite{BHR}):

\begin{thm} \label{taylor}
Given a cuspidal automorphic representation $\pi$ of ${\rm
GL}_2(\A_F)$ with $\pi_{\infty}$ isomorphic to the principal
series representation corresponding to $$\begin{pmatrix} t_1& *\\
0 & t_2
\end{pmatrix} \mapsto \left( \frac{t_1}{|t_1|} \right) \left(
\frac{|t_2|}{t_2} \right)$$  and cyclotomic central character
$\omega$ (i.e. $\omega^c=\omega$), let $\Sigma_{\pi}$ denote the set
of places above $p$, the primes where $\pi$ or $\pi^c$ is ramified,
and primes ramified in $F/\Q$.

Then there exists a continuous Galois representation
$$\rho_{\pi}:G_F \to {\rm GL}_2(\overline
F_{\mathfrak{p}})$$ such that if $v \notin \Sigma_{\pi}$, then
$\rho_{\pi}$ is unramified at $v$ and the characteristic
polynomial of $\rho_{\pi}(\mathrm{Frob}_v)$ is
$x^2-a_v(\pi)x+\omega(\mathfrak{P}_v)
\mathrm{Nm}_{F/\Q}(\mathfrak{P}_v),$ where $a_v(\pi)$ is the Hecke
eigenvalue corresponding to $T_v$. The image of the Galois
representation lies in $\mathrm{GL}_2(L)$ for a finite extension
$L$ of $F_{\mathfrak{p}}$ and the representation is absolutely
irreducible. \hspace{\fill} \qedsymbol
\end{thm}

\begin{rem}
\begin{enumerate}
\item Taylor relates $\pi$ to a low weight Siegel modular form via a theta lift
and uses the Galois representation attached to this form (via
pseudorepresentations and the Galois representations of
cohomological Siegel modular forms) to find $\rho_{\pi}$.
  \item Taylor had some additional technical assumption in \cite{T2} and
only showed the equality of Hecke and Frobenius polynomial outside
a set of places of zero density. For this strengthening of
Taylor's result see \cite{BHR}.
\item Conjecture 3.2 in \cite{CD} describes for cuspforms of general weight
a conjectural extension of Taylor's theorem.
\end{enumerate}
\end{rem}

Urban studied in \cite{U98} the case of ordinary automorphic
representations $\pi$, and together with results in \cite{U04} on
the Galois representations attached to ordinary Siegel modular
forms shows:
\begin{thm} [(Corollaire 2 of \cite{U04})] \label{Urban}
If $\pi$ is unramified at $\mathfrak{p}$ and  \textit{ordinary} at
$\mathfrak{p}$, i.e., $|a_{\mathfrak{p}}(\pi)|_p=1$, then the
Galois representation $\rho_{\pi}$ is ordinary at $\mathfrak{p}$,
i.e.,
$$\rho_{\pi}|_{G_{\mathfrak{p}}} \cong
\begin{pmatrix}\Psi_1& *\\ 0 & \Psi_2 \end{pmatrix},$$ where
$\Psi_2|_{I_{\mathfrak{p}}}=1$, and
$\Psi_1|_{I_{\mathfrak{p}}}=\mathrm{det}(\rho_{\pi})|_{I_{\mathfrak{p}}}
= \epsilon$.
\end{thm}

For $p$ inert we will need a stronger statement:
\begin{conjecture} \label{chrys}
If $\pi$ is unramified at $\mathfrak{p}$ then
$\rho_{\pi}|_{G_{\mathfrak{p}}}$ is crystalline.
\end{conjecture}
\noindent This conjecture extends Conjecture 3.2 in \cite{CD} and
would follow if one could prove the corresponding statement for
low weight Siegel modular forms.

\section{Selmer group and Eisenstein ideal} \label{SelEis}
In this section we define an Eisenstein ideal in a Hecke algebra
acting on cuspidal automorphic forms of ${\rm GL}_{2/F}$ and show
that its index gives a lower bound on the size of the Selmer group
of a Galois character.

 Let $\phi_1$ and $\phi_2$ be two Hecke
characters with infinity type $z$ and $z^{-1}$, respectively. Let
$\mathcal{R}$ be the ring of integers in the finite extension $L$ of
$F_{\mathfrak{p}}$ containing the values of the finite parts of
$\phi_i$ and $L^{\rm alg}(0,\phi_1/\phi_2)$. Denote its maximal
ideal by $\mathfrak{m}_{\mathcal{R}}$. Let $\Sigma_{\phi}$ be the
finite set of places dividing the conductors of the characters
$\phi_i$ and their complex conjugates and the places dividing
$pd_F$. Let $K_f=\prod_v K_v \subset G(\A_f)$ be a compact open
subgroup such that $K_v = {\rm GL}_2(\mathcal{O}_{v})$ if $v \notin
\Sigma_{\phi}$.

Because of the condition on the central character in Theorem
\ref{taylor} we assume that there exists a finite order Hecke
character $\eta$ unramified outside $\Sigma_{\phi}$ such that
\begin{equation} \label{eta} (\phi_1 \phi_2 \eta^2 )^c=\phi_1 \phi_2 \eta^2.
\end{equation}

Denote by $\mathbf{T}$ the $\mathcal{R}$-algebra generated by the
Hecke operators $T_v, v \notin \Sigma_{\phi}$ acting on
$S_2(K_f,\phi_1 \phi_2)$. Call the ideal $\mathbf{I}_{\phi}
\subseteq \mathbf{T}$ generated by
$$ \left\{ T_v-\phi_{1}(\mathfrak{P}_v){\rm Nm}(\mathfrak{P}_v) -
 \phi_{2}(\mathfrak{P}_v) | v \notin \Sigma_{\phi} \right\}
$$
the \textit{Eisenstein ideal associated to} $\phi=(\phi_1,
\phi_2)$.

Using the notation of \S \ref{s2.2}, we define Galois characters
\begin{eqnarray} \rho_1&=&\phi_{1,{\mathfrak{p}}} \epsilon, \nonumber\\
\rho_2&=&\phi_{2,{\mathfrak{p}}},\nonumber\\\rho&=&\rho_1 \otimes
\rho_2^{-1}. \nonumber
\end{eqnarray}
Notice that $\rho$ depends only on the quotient $\phi_1/\phi_2$.
Let $\Sigma_{\rho}$ be the set of places dividing $p$ and those
where $\rho$ is ramified. Our first main result is the following
inequality:
\begin{thm} \label{Selthm}
$$\mathrm{val}_p(\#{\rm Sel}^{\Sigma_{\phi}\backslash \Sigma_{\rho}}(F,\rho)) \geq
\mathrm{val}_p(\#(\mathbf{T}/\mathbf{I}_{\phi})).$$
\end{thm}

\begin{proof}
We can assume that
$$\mathbf{T}/\mathbf{I}_{\phi} \neq 0.$$
Let $\mathfrak{m}\subset \mathbf{T}$ be a maximal ideal containing
$\mathbf{I}_{\phi}$. Taking the completion with respect to
$\mathfrak{m}$ we write $$S_2(K_f,\phi_1
\phi_2)_{\mathfrak{m}}=\bigoplus_{i=1}^n V_{\pi_{i,f}}^{K_f},$$
where $V_{\pi_f}$ denotes the representation space of the (finite
part) of a cuspidal automorphic representation $\pi$.

By twisting the cuspforms by the finite order character $\eta$ of
(\ref{eta}) we can ensure that their central character is
cyclotomic. Hence we can apply Theorem \ref{taylor} to associate
Galois representations $\rho_{\pi_i \otimes \eta}:G_{\Sigma_{\phi}}
\to {\rm GL}_2(L_i)$ to each $\pi_i \otimes \eta, i=1 \ldots n$ for
some finite extensions $L_i/F_{\mathfrak{p}}$. Taking all of them
together (and untwisting by $\eta$) we obtain a continuous,
absolutely irreducible Galois representation
$$\rho_T:=\bigoplus_{i=1}^n \rho_{\pi_i \otimes \eta} \otimes \eta^{-1}: G_{\Sigma_{\phi}} \to {\rm
GL}_2(\mathbf{T}_{\mathfrak{m}} \otimes_{\mathcal{R}} L).$$ Here
we use that $\mathbf{T}_{\mathfrak{m}} \otimes_{\mathcal{R}}
L=\prod_{i=1}^n L_i$, which follows from the strong multiplicity
one theorem. We have an embedding
$$\mathbf{T}_{\mathfrak{m}} \hookrightarrow \prod_{i=1}^n L_i
$$
$$T_v \mapsto ((a_v(\pi_i)),$$ where $a_v(\pi_i)$ is the
$T_v$-eigenvalue of $\pi_i$. The coefficients of the
characteristic polynomial ${\rm char}(\rho_T)$ lie in
$\mathbf{T}_{\mathfrak{m}}$ and by the Chebotarev density theorem
$${\rm char}(\rho_T) \equiv {\rm char}(\rho_1 \oplus \rho_2) \mod{\mathbf{I}_{\phi}}.$$

For any finite free $\Tm \otimes L$-module $\mathcal{M}$, any
$\Tm$-submodule $\mathcal{L} \subset \mathcal{M}$ that is finite
over $\Tm$ and such that $\mathcal{L} \otimes L=\mathcal{M}$ is
called a \textit{$\Tm$-lattice}. Specializing to our situation
Theorem 1.1 of \cite{U01} we get:
\begin{thm} [(Urban)] \label{ThmUrban}
Given a Galois representation $\rho_T$ as above there exists a
$G_{\Sigma_{\phi}}$-stable $\Tm$-lattice $\mathcal{L} \subset (\Tm
\otimes L)^2$  such that $G_{\Sigma_{\phi}}$ acts on
$\mathcal{L}/\mathbf{I}_{\phi} \mathcal{L}$ via the short exact
sequence
$$0 \to \mathcal{R}_{\rho_1} \otimes_{\mathcal{R}} (N/\mathbf{I}_{\phi}) \to
\mathcal{L}/\mathbf{I}_{\phi} \mathcal{L} \to \mathcal{R}_{\rho_2}
\otimes_{\mathcal{R}} (\Tm/\mathbf{I}_{\phi}) \to 0,$$ where $N
\subset \Tm \otimes L$ is a $\Tm$-lattice with ${\rm val}_p(\#
\Tm/\mathbf{I}_{\phi}) \leq {\rm val}_p(\# N/\mathbf{I}_{\phi}N) <
\infty$ and no quotient of $\mathcal{L}$ is isomorphic to
$\overline{\rho}_1:=\rho_1 \mod{\mathfrak{m}_{\mathcal{R}}}$.
\end{thm}
\begin{proof}
Note that the  $\mathcal{R}$-algebra map $\mathcal{R}
\twoheadrightarrow \Tm/\mathbf{I}_{\phi}$ is surjective and that
$\mathcal{L}/I_{\phi} \mathcal{L} \cong \mathcal{L}
\otimes_{\mathcal{R}} \Tm/\mathbf{I}_{\phi}$. Hence this short exact
sequence recovers the one in the statement of Theorem 1.1 of
\cite{U01}. For the statement about ${\rm val}_p(\#
N/\mathbf{I}_{\phi}N)$ see \cite{U01} p. 519 and use that any
$\mathcal{R}$-submodule of $\Tm/\mathbf{I}_{\phi}$ or
$N/\mathbf{I}_{\phi}N$ is a $\Tm$-submodule.

See \cite{Be} \S 7.3.2 for an alternative construction of such a
lattice using arguments of Wiles (\cite{Wi86} and \cite{Wi}).
\end{proof}

Using the properties of the Galois representations attached to
cuspforms listed in \S \ref{Galreps} we can now conclude the proof
of Theorem \ref{Selthm} by similar arguments as in \cite{S04} and
\cite{U01}. To ease notation we put
$\mathcal{T}:=N/\mathbf{I}_{\phi}$ and $\Sigma:=\Sigma_{\phi}$.

Identifying $\mathcal{R}_{\rho}$ with ${\rm
Hom}_{\mathcal{R}}(\mathcal{R}_{\rho_2},\mathcal{R}_{\rho_1})$ and
writing $s:\mathcal{R}_{\rho_2} \otimes \Tm/\mathbf{I}_{\phi} \to
\mathcal{L} \otimes \Tm/\mathbf{I}_{\phi}$ for the section as
$\Tm/\mathbf{I}_{\phi}$-modules we define a 1-cocycle
$c:G_{\Sigma} \to \mathcal{R}_{\rho} \otimes \mathcal{T}$ by
$$c(g)(m)= \text{ the image of } s(m)-g.s(g^{-1}.m) \text{ in }
\mathcal{R}_{\rho_1}\otimes \mathcal{T}.$$

Consider the $\mathcal{R}$-homomorphism $$\varphi:{\rm
Hom}_{\mathcal{R}}(\mathcal{T}, \mathcal{R}^{\vee}) \to
H^1(G_{\Sigma}, W_{\rho}), \; \varphi(f)= \text{ the class of } (1
\otimes f) \circ c.$$ We will show that

\begin{enumerate}[(i)]
  \item ${\rm im}(\varphi) \subset {\rm Sel}^{\Sigma \backslash \Sigma_{\rho}}(F, \rho)$,
  \item ${\rm ker}(\varphi)=0$.
\end{enumerate}
From (i) it follows that $${\rm val}_p(\#{\rm Sel}^{\Sigma
\backslash \Sigma_{\rho}}(F, \rho)) \geq {\rm val}_p(\# {\rm
im}(\varphi)).$$ From (ii) it follows that
\begin{eqnarray*}{\rm val}_p(\# {\rm im}(\varphi)) &\geq&
{\rm val}_p(\# {\rm Hom}_{\mathcal{R}}(\mathcal{T},
\mathcal{R}^{\vee}))\\ &=&{\rm val}_p(\# \mathcal{T})
\\ &\geq& {\rm val}_p(\# \Tm/\mathbf{I}_{\phi}).\end{eqnarray*}

 For (i) we have to show that the
 conditions of the Selmer group at $v \mid p$ are satisfied:
 For split $p$
it suffices to prove that the extension in Theorem \ref{ThmUrban}
is
 split when considered as an extension of $\Tm[G_{\mathfrak{p}}]$-modules, because then the class
 in $H^1(G_{\mathfrak{p}}, \mathcal{R}_{\rho} \otimes \mathcal{T})$ determined by $c$ is the zero class.
In this case the Hecke eigenvalues $a_{\mathfrak{p}}(\pi_i) \equiv
p \cdot \phi_1(\mathfrak{p})+ \phi_2(\mathfrak{p}) \not \equiv 0
\mod{\mathfrak{m}_R}$, hence the cuspforms $\pi_i \otimes \eta$
are ordinary at $\mathfrak{p}$, so Theorem \ref{Urban} applies and
$\rho_T$ is ordinary. Observing that the Hodge-Tate weights at
$\mathfrak{p}$ of $\rho_1$ and $\rho_2$ are  0 and 1,
respectively, the splitting of the extension  as
$\Tm[G_{\mathfrak{p}}]$-modules follows from comparing the basis
given by Theorem \ref{ThmUrban} with the one coming from
ordinarity.

For inert $p$ we observe that by Conjecture \ref{chrys} the
$\rho_{\pi_i}$ are all crystalline, which implies that the class
determined by $c$ is crystalline.

To prove (ii) we first observe that for any $f \in {\rm
Hom}_{\mathcal{R}}(\mathcal{T}, \mathcal{R}^{\vee})$, ${\rm
ker}(f)$ has finite index in $\mathcal{T}$ since $\mathcal{T}$ is
a finite $\mathcal{R}$-module and so $f \in {\rm
Hom}_{\mathcal{R}}(\mathcal{T},
\mathcal{R}^{\vee}[\mathfrak{m}_{\mathcal{R}}^n])$ for some $n$.
Suppose now that $f \in {\rm ker}(\varphi)$. We claim that the
class of $c$ in $H^1(G_{\Sigma}, \mathcal{R}_{\rho}
\otimes_{\mathcal{R}} \mathcal{T}/{\rm ker}(f))$ is zero. To see
this, let $X=\mathcal{R}^{\vee}/{\rm im}(f)$ and observe that
there is an exact sequence
$$ H^0(G_{\Sigma}, \mathcal{R}_{\rho} \otimes_{\mathcal{R}} X) \to H^1(G_{\Sigma}, \mathcal{R}_{\rho} \otimes_{\mathcal{R}} \mathcal{T}/{\rm
ker}(f)) \to H^1(G_{\Sigma}, \mathcal{R}_{\rho}
\otimes_{\mathcal{R}} \mathcal{R}^{\vee}).$$ Since $f \in {\rm
ker}(\varphi)$ and the second arrow in the sequence comes from the
inclusion $\mathcal{T}/{\rm ker}(\varphi) \hookrightarrow
\mathcal{R}^{\vee}$ induced by $f$, the image in the right module
of the class of $c$ in the middle is zero. Our claim follows
therefore if the module on the left is trivial. But the dual of
this module is a subquotient of ${\rm
Hom}_{\mathcal{R}}(\mathcal{R}_{\rho}, \mathcal{R})$ on which
$G_{\Sigma}$ acts trivially. By assumption, however,
$\mathcal{R}_{\rho}$ is a rank one module on which $G_{\Sigma}$
acts non-trivially.

Suppose in addition that $f$ is non-trivial, i.e., ${\rm ker}(f)
\subsetneq \mathcal{T}$. Note that any $\mathcal{R}$-submodule of
$\mathcal{T}$ is actually a $\Tm$-submodule since $\mathcal{R}
\twoheadrightarrow \Tm/\mathbf{I}_{\phi}$. Hence there exists a
$\Tm$-module $A$ with ${\rm ker}(f) \subset A \subset \mathcal{T}$
such that $\mathcal{T}/A \cong
\mathcal{R}/\mathfrak{m}_{\mathcal{R}}$. From our claim it follows
that the $\Tm[G_{\Sigma}]$-extension
$$0 \to \mathcal{R}_{\rho_1} \otimes_{\mathcal{R}}
\mathcal{R}/\mathfrak{m}_{\mathcal{R}} \cong \mathcal{R}_{\rho_1}
\otimes_{\mathcal{R}} \mathcal{T}/A \to
\mathcal{L}/(\mathcal{R}_{\rho_1} \otimes_{\mathcal{R}} A) \to
\mathcal{L}/\mathcal{L}_1 \to 0$$ is split. But this would give a
$\Tm[G_{\Sigma}]$-quotient of $\mathcal{L}$ isomorphic to
$\overline{\rho}_{1}$, which contradicts one of the properties of
the lattice constructed by Urban. Hence ${\rm ker}(\varphi)$ is
trivial.
\end{proof}

The following Lemma will later provide us with the finite order
character $\eta$ of (\ref{eta}) used in the twisting above.
\begin{lem}\label{anticycchar}
If $\chi=\phi_1/\phi_2$ satisfies $\chi^c=\overline{\chi}$ then
there exists a finite order character $\eta$ unramified outside
$\Sigma_{\phi}$ such that $(\phi_1 \phi_2 \eta^2 )^c=\phi_1 \phi_2
\eta^2$.
\end{lem}

\begin{proof}
We claim that there exists a Hecke character $\mu$ unramified
outside $\Sigma_{\phi}$ such that $$\chi=\mu \overline \mu^c.$$
Given such a character $\mu$ we then define $\eta=(\mu
\phi_2)^{-1}$.

In the Lemma on p.81 of \cite{G83} Greenberg defines a Hecke
character $\mu_G:F^* \backslash \A_F^* \to \C^*$ of infinity type
$z^{-1}$ such that $\mu_G^c=\overline \mu_G$ and $\mu_G$ is
ramified exactly at the primes ramified in $F/\Q$. It therefore
suffices to prove the claim for the finite order character
$$\chi':=\chi \mu_G^{2}=\chi \mu_G (\overline \mu_G^c).$$

By assumption we have that
$$\chi' \equiv 1 \text{ on } \mathrm{Nm}_{F/\Q}(\A_F^*) \subset \A^*_{\Q}
\subset \A_F^*.$$ Thus $\chi'$ restricted to $\Q^* \backslash
\A_{\Q}^*$ is either the quadratic character of $F/\Q$ or trivial.
Since our finite order character has trivial infinite component,
$\chi'$ has to be trivial on $\Q^* \backslash  \A_{\Q}^*$.
Hilbert's Theorem 90 then implies that there exists $\mu$ such
that $\chi'=\mu/\mu^c$.

To control the ramification we analyze this last step closer:
$\chi'$ factors through $\A_F^* \to A$, where $A$ is the subset of
$\A_F^*$ of elements of the form $x/x^c$ and the map is $x \mapsto
x/x^c$. If $y \in A \cap F^*$ then $y$ has trivial norm and so by
Hilbert's Theorem 90, $y=x/x^c$ for some $x \in F^*$. Thus the
induced character $A \to \C^*$ vanishes on $A \cap F^*$. This
implies that there is a continuous finite order character $\mu:
F^* \backslash \A_F^* \to \C^*$ which restricts to this character
on $A$ and $\chi'=\mu/\mu^c$ (this argument is taken from the
proof of Lemma 1 in \cite{T2}).

By the following argument we can further conclude that the induced
character vanishes on $A \cap \prod_{v \notin \Sigma_{\phi}}
\mathcal{O}_v^*$ and therefore find $\mu$ on $F^* \backslash
\A_F^*/ \C^* \prod_{v \notin \Sigma_{\phi}} \mathcal{O}_v^*$
restricting to the character $A \to \C^*$: Writing $U_{F,{\ell}}=
\prod_{v |{\ell}} \mathcal{O}_v^*$ for a prime $\ell$ in $\Q$ we
have
$$H^1({\rm Gal}(F/\Q), \prod_{v \notin \Sigma_{\phi}}
\mathcal{O}_v^*) \cong \prod_{\ell \notin \Sigma_{\phi}} H^1({\rm
Gal}(F/\Q), U_{F,{\ell}}),$$ where ``$\ell \notin \Sigma_{\phi}$"
denotes those $\ell \in \Z$ such that $v \mid \ell \Rightarrow v
\notin  \Sigma_{\phi}$. For the isomorphism we use that $$v \in
\Sigma_{\phi} \Rightarrow \overline v \in \Sigma_{\phi}.$$ In
fact, all these groups are trivial since all $\ell \notin
\Sigma_{\phi}$ are unramified in $F/\Q$ and so
$$H^1({\rm Gal}(F/\Q), U_{F,{\ell}}) \cong H^1(G_v,
\mathcal{O}_v^*)=1.$$ If $y \in A \cap \prod_{v \notin
\Sigma_{\phi}} \mathcal{O}_v^*$ then $y$ has trivial norm in
$\prod_{v\notin \Sigma_{\phi}} \mathcal{O}_v^*$. But as shown, its
first Galois cohomology group is trivial so there exists $x \in
\prod_{v \notin \Sigma_{\phi}} \mathcal{O}_v^* \cap \A_F^*$ such
that $y=x/x^c$. Since $\chi'$ is unramified outside
$\Sigma_{\phi}$ the image of $y$ under the induced character
therefore equals $\chi'(x)=1$, as claimed above.
\end{proof}

\section{Cohomological congruences} \label{s5}

In \cite{Be06} we constructed a class $\Eis$ in the cohomology of
a symmetric space associated to ${\rm GL}_{2/F}$ with integral
non-zero restriction to the boundary of the Borel-Serre
compactification that is an eigenvector for the Hecke operators at
almost all places. By a result of Harder one knows that $\Eis$ is
rational. The main result of \cite{Be06} is a lower bound on its
denominator (defined in  (\ref{denominator}) below) in terms of
the $L$-value of a Hecke character. In this section we show that
if there exists an integral cohomology class with the same
restriction to the boundary as $\Eis$ then there exists a
congruence modulo the $L$-value between $\Eis$, multiplied by its
denominator, and a cuspidal cohomology class. This implies that
the $L$-value bounds the index of the Eisenstein ideal from below.
\subsection{The Eisenstein cohomology set-up}
Assume from now on in addition that $p>3$. Recall the notation and
definitions introduced in Section \ref{s2.5}. Let $\Ox$ denote the
ring of integers in the finite extension $\Fx$ of $F_{\mathfrak{p}}$
obtained by adjoining the values of the finite part of both $\phi_i$
and $L^{\rm alg}(0,\phi_1/\phi_2)$. We write
$$\widetilde H^1(X, \mathcal{R}):=H^1(X, \mathcal{R})_{\rm
free}=\mathrm{im} (H^1(X,\mathcal{R}) \to H^1(X,L))$$ for
$X=S_{K_f}$ or $\partial \overline S_{K_f}$. For $c \in
H^1(S_{K_f},\Fx)$ define the denominator (ideal) by \begin{equation}
\label{denominator} \delta(c):=\{a \in \Ox: a \cdot c \in \widetilde
H^1(S_{K_f},\Ox) \}.\end{equation} Recall the long exact sequence
$$\ldots \to H^1_c(S_{K_f},R) \to H^1(S_{K_f},R)
\overset{\rm{res}} \to H^1(\partial \overline S_{K_f},R) \to
H^2_c(S_{K_f},R) \to \ldots$$ for any $\Ox$-algebra $R$.

\noindent \textit{The set-up.} Suppose we are given a pair of
Hecke characters $\phi=(\phi_1, \phi_2)$ and a class $\Eis \in
H^1(S_{K_f},\Fx)$ satisfying the following properties:
\begin{enumerate} [(E1)]
  \item  The image of $\Eis$ under ${\rm res}$ lies in
  $\widetilde H^1(\partial \overline S_{K_f},\Ox)$.
  \item For all places $v$ outside the conductors of the $\phi_i$
  the class $\Eis$ is an eigenvector for the Hecke operator $T_{\pi_v}=[K_f \begin{pmatrix} \pi_v &0\\0&1
\end{pmatrix} K_f]$ with eigenvalue $$\phi_{2}(\mathfrak{P}_v)+ {\rm Nm}(\mathfrak{P}_v)
\phi_{1}(\mathfrak{P}_v).$$
    \item The central character of $\Eis$
is given by $\phi_1 \phi_2$, i.e., the Hecke operators
$$S_{\pi_v}=[K_f \begin{pmatrix} \pi_v &0\\0&\pi_v
\end{pmatrix} K_f]$$ act on it by multiplication by $(\phi_1
\phi_2)(\pi_v)$.
\item  The denominator of $\Eis$ is bounded below by $L^{\mathrm{alg}}(0,
\phi_1/\phi_2)$, i.e., $$\delta(\Eis) \subseteq
(L^{\mathrm{alg}}(0, \phi_1/\phi_2)).$$
\end{enumerate}
Suppose we are also given:
\begin{enumerate} [(H1)]
  \item There exists $c_{\phi}
  \in \widetilde H^1(S_{K_f},\Ox)$ with $${\rm res}(c_{\phi})={\rm
  res}(\Eis) \in \widetilde H^1(\partial \overline S_{K_f}, \Ox).$$
  \item There exists an idempotent $e_{\omega}$ acting on
  $H^1(S_{K_f},\C)$ such that $S_{\pi_v} e_{\omega}= (\phi_1
  \phi_2)(\pi_v) e_{\omega}$ for $v$ not dividing the conductors of the $\phi_i$.
\end{enumerate}
\vspace{0.5cm}

The following provides a bound on the congruence module introduced
in the previous section:
\begin{prop} \label{Eisideal}
Given the above setup there is an $\Ox$-algebra surjection
$$
\mathbf{T}/\mathbf{I}_{\phi} \twoheadrightarrow \Ox/ \left(
L^{\mathrm{int}}(0,\phi_1/\phi_2)  \right).
$$
\end{prop}

\begin{proof}
Put $\widetilde H^1_!(S_{K_f}, \mathcal{R})=H^1_!(S_{K_f}, L) \cap
\widetilde H^1(S_{K_f}, \mathcal{R})$ and $\omega=\phi_1 \phi_2$.
Under the Eichler-Shimura-Harder isomorphism (see (\ref{ESH})) we
have
$$ e_{\omega} H^1_!( S_{K_f}, \C) \cong S_2(K_f, \omega).$$
Hence the Hecke algebra $\mathbf{T}$ from Section \ref{SelEis} is
isomorphic to
 the $\Ox$-subalgebra of $$
{\rm End}_{\Ox}(e_{\omega}  \widetilde H^1_!( S_{K_f} , \Ox))$$
generated by the Hecke operators $T_{\pi_v}$ for all primes $v
\notin \Sigma_{\phi}$ and we will identify the two.

Note that for $c_{\phi} \in \widetilde H^1(S_{K_f},\Ox)$ given by
(H1) we have
$${\rm res}(e_{\omega} c_{\phi})=e_{\omega} {\rm res}(c_{\phi})=e_{\omega}{\rm
res}(\Eis)={\rm res}(\Eis)$$ since $S_v (\Eis)= \omega(\pi_v)\Eis$
by (E3).

Without loss of generality, we can assume that $\delta(\Eis) \neq
\Ox$; there is nothing to prove otherwise by (E4). Let $\delta$ be
a generator of $\delta(\Eis)$. Then $\delta \cdot \Eis$ is an
element of an $\Ox$-basis of $e_{\omega} \widetilde
H^1(S_{K_f},\Ox)$. By construction, $c_0:= \delta \cdot
(e_{\omega} c_{\phi} - \Eis) \in e_{\omega} H^1_!(S_{K_f},L)$ is a
nontrivial element of an $\Ox$-basis of $e_{\omega} \widetilde
H^1_!(S_{K_f},\Ox)$. Extend $c_0$ to an $\Ox$-basis $c_0, c_1,
\ldots c_d$ of $e_{\omega} \widetilde H^1_!(S_{K_f},\Ox)$. For
each $t\in \mathbf{T}$ write
$$t(c_0)= \sum_{i=0}^{d} a_i(t) c_i, \, a_i(t) \in \Ox.$$
Then \begin{equation} \label{map} \mathbf{T} \to \Ox/(\delta), \;
t \mapsto a_0(t) \mod {\delta} \end{equation} is an $\Ox$-module
surjection. We claim that it is independent of the $\Ox$-basis
chosen and that it is a homomorphism of $\Ox$-algebras with the
Eisenstein ideal $\mathbf{I}_{\phi}$ contained in its kernel. To
prove this it suffices to check that each $a_0(T_{\pi_v} -
\phi_{2}(\mathfrak{P}_v) - {\rm{Nm}}(\mathfrak{P}_v)
\phi_{1}(\mathfrak{P}_v))$, $v \notin \Sigma_{\phi}$ is contained
in $\delta \Ox$. This is an easy calculation using (E2). Since
$\Ox/(\delta) \twoheadrightarrow \Ox/(L^{\mathrm{int}}(0,\chi))$
by (E4), this concludes the proof.
\end{proof}

In the following sections, we will indicate how to produce the
elements in the set-up of the Proposition. Under certain conditions
on the character $\phi_i$, to be reviewed in Section
\ref{Eisreview}, we constructed in \cite{Be06} a class $\Eis$
satisfying (E1)-(E4) using Harder's Eisenstein cohomology.
Assumption (H2) is of a technical nature and will be discussed in
Section \ref{idempotent}. We are interested in controlling the
central character via (H2) because of the restriction in Theorem
\ref{taylor}. The most difficult ingredient to procure is (H1), see
Sections \ref{tproblem} and \ref{s6}.

\begin{rem}
\begin{enumerate}
  \item As already remarked in the introduction, the
  constant cohomology coefficients above can be replaced by  coefficient systems
  arising from finite dimensional representations of ${\rm GL}_{2/F}$.
  See \cite{Be06} \S 2.4 and 3.1 for the
  necessary modifications.
  \item Note also that except for the explicit Hecke operators we did not use any
information specific to ${\rm GL}_{2/F}$, i.e., $S_{K_f}$ could be
replaced by a symmetric space associated to a different group $G$
and $\phi$ by a tuple of automorphic forms on the Levi part of a
parabolic subgroup of $G$. Since the analytic theory of Eisenstein
cohomology has been developed for a wide variety of groups, and
rationality results are known, e.g. for ${\rm GL}_n$ by the work
of \cite{FS}, we hope that our techniques generalize to these
groups.
\end{enumerate}
\end{rem}

\subsection{Construction of Eisenstein class} \label{Eisreview}
Following Harder we constructed in \cite{Be06} Eisenstein
cohomology classes in the Betti cohomology group
$H^1(S_{K_f},\C)$. Given a pair of Hecke characters
$\phi=(\phi_1,\phi_2)$ with $\phi_{1,\infty}(z)=z$ and
$\phi_{2,\infty}(z)=z^{-1}$ these depend on a choice of a function
$\Psi_{\phi_f}$ in the induced representation
$$V_{\phi_f, \C}^{K_f}=\{\Psi:G(\A_f) \to \C | \Psi(bg)=\phi_f(b) \Psi(g)
\forall b \in B(\A_f), \Psi(gk)=\Psi(g) \forall k \in K_f \}.$$ In
the notation of \cite{Be06} we take $K_f=K_f^S$ and
$\Psi_{\phi_f}=\Psi^0_{\phi}$. We recall the definition of the
compact open $K_f$: Denote by $S$ the finite set of places where
both $\phi_i$ are ramified, but $\phi_1/\phi_2$ is unramified.
Write $\mathfrak{M}_i$ for the conductor of $\phi_i$. For an ideal
$\mathfrak{N}$ in $\mathcal{O}$ and a finite place $v$ of $F$ put
$\mathfrak{N}_v=\mathfrak{N} \mathcal{O}_v$. We define
$$K^1(\mathfrak{N})=\left \{
\begin{pmatrix}
  a & b \\
  c & d
\end{pmatrix}
\in {\rm GL}_2(\widehat{\mathcal{O}}), a-1, c \equiv 0 \,
\mod{\mathfrak{N}}  \right \},$$
$$K^1(\mathfrak{N}_v)=\left \{
\begin{pmatrix}
  a & b \\
  c & d
\end{pmatrix}
\in {\rm GL}_2(\mathcal{O}_v), a-1, c \equiv 0 \,
\mod{\mathfrak{N}_v} \right \},$$ and $$U^1(\mathfrak{N}_v)=\{k
\in {\rm GL}_2(\mathcal{O}_v): {\rm det}(k) \equiv 1
\mod{\mathfrak{N}_v} \}.$$ Now put
$$K_f:=\prod_{v \in S} U^1(\mathfrak{M}_{1,v}) \prod_{v \notin S}
K^1((\mathfrak{M}_1 \mathfrak{M}_2)_v).$$ The exact definition of
$\Psi_{\phi_f}$ will not be required in the following; we refer
the interested reader to \cite{Be06} Section 3.2. Let $\Eis$ be
the cohomology class denoted by $[{\rm
Eis}(\Psi^0_{(\phi_1,\phi_2)_f})]$ in \cite{Be06}.

The rationality of $\Eis$, i.e., the fact that $\Eis \in
H^1(S_{K_f},\Fx)$ was proven by Harder, see  \cite{Be06} Proposition
13. Properties (E2) and (E3) are satisfied by construction, see
\cite{Be06} Lemma 9. The integrality of the constant term (E1) is
analyzed in \cite{Be06} Proposition 16. The main result of
\cite{Be06} is the bound on the denominator (E4). The latter two
require certain conditions on the characters $\phi_1$ and $\phi_2$.
However, since in the combination of Theorem \ref{Selthm} and
Proposition \ref{Eisideal} the main object of interest is the
character $\chi=\phi_1/\phi_2$, we will from now on focus on $\chi$
and view $\phi_1$ and $\phi_2$ as auxiliary.

\begin{thm}[(\cite{Be06} Proposition 16, Theorem 29)] \label{thmdenom}
Let $\chi$ be a Hecke character $\chi$ of infinity type $z^2$ with
conductor $\mathfrak{M}$ coprime to $(p)$. Assume in addition that
either

(i) $p$ splits in $F$, $\chi$ has split conductor, and
$\frac{L(0,\overline \chi)}{L(0,\chi)} \in \Ox$

or

 (ii) $\chi^c=\overline \chi$, no
ramified primes  divide $\mathfrak{M}$ and no inert primes
congruent to $-1 \mod{p}$ divide $\mathfrak{M}$ with multiplicity
one, and
$$\omega_{F/\Q}(\mathfrak{M}) \frac{\tau(\tilde \chi)}{\sqrt{{\rm Nm}(\mathfrak{M})}}=1,$$
where $\omega_{F/\Q}$ is the quadratic Hecke character associated
to the extension $F/\Q$ and $\tau(\tilde \chi)$ the Gauss sum of
the unitary character $\tilde \chi:=\chi/|\chi|$.

Then there exists a factorization $\chi=\phi_1/\phi_2$ such that
$\Eis$ satisfies (E1)-(E4). \hfill \qedsymbol
\end{thm}

\begin{rem}
\begin{enumerate}
  \item Proposition 16 of \cite{Be06} shows that $\frac{L(0,\overline \chi)}{L(0,\chi)} \in
  L$.
  \item By considering non-constant coefficient systems
  \cite{Be06},
  in fact, proves this for characters $\chi$ of infinity type
  $z^{m+2} \overline z^{-m}$ for $m \in \N_{\geq 0}$, if $p$ splits
  in $F/\Q$.
\end{enumerate}
\end{rem}

\subsection{Existence of idempotent (H2)} \label{idempotent}
\begin{lem} \label{H2}
Let $K_f$ be the compact open defined in Section \ref{Eisreview}.
If $p \nmid \#{\rm Cl}_{\mathfrak{M} \mathfrak{M}_1}(F)$ then (H2)
is satisfied.
\end{lem}

\begin{proof}
By \cite{U98} \S 1.2 and \S 1.4.5 the action of the diamond
operators $S_{\pi_v} , v \nmid {\rm cond}(\phi_i)$ on $H^1(S_{K_f},
\C)$ is determined by the class in ${\rm Cl}_{\mathfrak{M}
\mathfrak{M}_1}(F)$ of the ideal determined by $\pi_v$ and induces
an $\Ox$-linear action of ${\rm Cl}_{\mathfrak{M}
\mathfrak{M}_1}(F)$ on $H^1(S_{K_f}, \C)$. Here we use that $$K_f
\supset K(\mathfrak{M}
\mathfrak{M}_1):=\left\{\begin{pmatrix}a&b\\c&d
\end{pmatrix} \in {\rm GL}_2(\hat{\mathcal{O}}): \begin{pmatrix}a&b\\c&d
\end{pmatrix} \equiv \begin{pmatrix}1&0\\0&1
\end{pmatrix} \mod{\mathfrak{M} \mathfrak{M}_1} \right\}.$$
By assumption the ray class group has order prime to $p$, so
$\Ox[{\rm Cl}_{\mathfrak{M} \mathfrak{M}_1}(F)]$ is semisimple.
For $\omega:=\phi_1 \phi_2$, which can be viewed as a character of
${\rm Cl}_{\mathfrak{M} \mathfrak{M}_1}(F)$, let $e_{\omega}$ be
the idempotent associated to $\omega$, so that $S_v e_{\omega} =
\omega(\pi_v) e_{\omega}$.
\end{proof}

\begin{rem}
By enlarging $K_f$ the condition $p \nmid \#
(\mathcal{O}/\mathfrak{M}\mathfrak{M}_1)^*$ can be weakened to the
order of $\phi_i|_{\hat{\mathcal{O}}^*}$ being coprime to $p$, see
\cite{Be} \S 6.1.
\end{rem}

\subsection{Torsion problem (H1)} \label{tproblem}
Hypothesis (H1) is related to the question of the occurrence of
torsion classes in $H^2_c(S_{K_f},\Ox)$. This problem does not arise
for ${\rm GL}_{2/\Q}$ because no such classes exist with the Hecke
eigenvalues under consideration (see \cite{CS} Lemma 6.1,
\cite{HaPi}). By Lefschetz duality (see \cite{G67} (28.18) or
\cite{Mau} Theorem 5.4.13)
$$H^2_c(S_{K_f},\Ox) \cong H_1(S_{K_f},\Ox),$$ so this question
reduces to the problem of torsion in $\Gamma^{\rm ab}$ for
arithmetic subgroups $\Gamma \subset G(\Q)$. This has been studied
in \cite{EGM82}, \cite{SV}, and \cite{GSch} (see also \cite{EGM}
\S 7.5). An arithmetic interpretation or explanation for the
torsion has not been found yet in general (but see \cite{EGM82}
for examples in the case of $\Q(\sqrt{-1}$)). Based on computer
calculations \cite{GSch} (2) suggests that for $\Gamma \subset
{\rm PSL}_2(\mathcal{O})$ apart from $2$ and $3$ only primes less
than or equal to $\frac{1}{2}[{\rm PSL}_2(\mathcal{O}):\Gamma]$
occur in the torsion of $\Gamma^{\rm ab}$. Even restricting to the
ordinary part there can be torsion, see \cite{Ta}\S 4. In all
cases calculated so far, ${\rm PSL}_2(\mathcal{O})^{\rm ab}$ has
only $2$ or $3$-torsion (see also \cite{Swan}, \cite{Berk}) but
this is not known in general, hence our different approach in the
following section, where we will prove:
  \begin{prop} \label{thmunram}
Let $\chi=\phi_1/\phi_2$ be an unramified Hecke character of
infinity type $z^2$. Assume that $1$ is the only unit in
$\mathcal{O}^*$ congruent to $1$ modulo the conductor of $\phi_1$.
If (E1) holds for $K_f$ and $\Eis$ as defined in Section
\ref{Eisreview}, then (H1) is satisfied.
  \end{prop}

\subsection{Congruence results} \label{s5.5}
We will summarize in this section the conditions under which we
can procure the ingredients for Proposition \ref{Eisideal} and
hence prove the existence of cohomological congruences.

\begin{thm} \label{condthm}
Assume $p$ splits in $F/\Q$. Let $\chi$ be a Hecke character $\chi$
of infinity type $z^2$ with split conductor $\mathfrak{M}$ coprime
to $(p)$. Assume $\frac{L(0,\overline \chi)}{L(0,\chi)} \in \Ox$ and
$$p \nmid \#(\mathcal{O}/\mathfrak{M})^* \cdot \#{\rm Cl}(F).$$ Let
$q> \#(\mathcal{O}^*)$  be any rational prime coprime to
$(p)\frak{M}$ and split in $F$ such that $p \nmid q-1$ and
$\frak{q}$ a prime of $F$ dividing $(q)$. If $H^2_c(S_{\tilde
K_f},\Z_p)_{\rm torsion}=0,$ where $$\tilde K_f=\{k \in
K^1(\mathfrak{M}): {\rm det}(k) \equiv 1 \mod{\frak{q}} \},$$ then
there exists a pair of characters $\phi=(\phi_1,\phi_2)$ with
$\chi=\phi_1/\phi_2$ such that  there is an $\Ox$-algebra
surjection
$$
\mathbf{T}/\mathbf{I}_{\phi} \twoheadrightarrow \Ox/ \left(
L^{\mathrm{int}}(0,\chi)  \right).
$$
\end{thm}

\begin{rem} As noted before, this
result is true, in fact, for characters $\chi$ of infinity type
$z^{m+2} \overline z^{-m}$ for $m \in \N_{\geq 0}$.
\end{rem}

\begin{proof}
Let $\phi_1$ be a  Hecke character with conductor $\frak{q}$ of
infinity type $z$ (for existence see, e.g. Lemma 24 of \cite{Be06}).
This is the character used in the proof of \cite{Be06} Theorem 29
and $K_f=\widetilde K_f$ for this pair $(\phi_1,\phi_1/\chi)$, so
the theorem follows from Proposition \ref{Eisideal}, Theorem
\ref{thmdenom}, and Lemma \ref{H2}.
\end{proof}

A similar result can be deduced for characters $\chi$ satisfying
$\overline \chi=\chi^c$ by taking as $\phi_1$ the character used
in the proof of \cite{Be06} Theorem 29. Its construction is more
involved and we refer the reader to the account in \cite{Be06} \S
5.2 and 5.3. The conductor $\mathfrak{M}_1$ of $\phi_1$ in this
case is given by $r \mathcal{D}$ for $\mathcal{D}$ the different
of $F$ and $r \in \Z$ any integer coprime to $(p) \mathfrak{M}$,
but such that no inert prime congruent to $-1$ modulo $p$ divides
$r$ with multiplicity one.

To be able to apply Lemma \ref{tameinertia} in Section \ref{s7} and
to satisfy the assumption in Proposition \ref{thmunram} we want to
impose the following extra condition on the conductor
$\mathfrak{M}_1$:
\begin{enumerate} [($\phi$)]
  \item $1$ is the only unit in
$\mathcal{O}^*$ congruent to $1$ modulo $\mathfrak{M}_1$, and
$$v \mid
  \mathfrak{M}_1 \Rightarrow v=\overline v \text{ and }
  \#\mathcal{O}_v/\mathfrak{P}_v \not \equiv \pm 1 \mod{p},$$
\end{enumerate}

We therefore assume in addition that $p \nmid \#{\rm Cl}(F)$ and
that $\ell \not \equiv \pm 1 \mod{p}$ for $\ell \mid d_F$. Also we
choose $r$ appropriately such that $p \nmid (\mathcal{O}/r)^*$ and
that $(\phi)$ holds.

We leave the counterpart of Theorem \ref{condthm} for characters
$\chi$ satisfying $\overline \chi=\chi^c$ to the interested reader
and instead give the following result, which does not require
torsion freeness. By Lemma \ref{unanticyc} unramified characters
$\chi$ satisfy $\chi^c=\overline \chi$, so we deduce from
Proposition \ref{Eisideal} together with Lemma \ref{H2}, Theorem
\ref{thmdenom}, and Proposition \ref{thmunram}:
\begin{thm} \label{unramcor}
Assume in addition that $p \nmid \#{\rm Cl}(F)$ and that $\ell
\not \equiv \pm 1 \mod{p}$ for $\ell \mid d_F$. Let $\chi$ be an
unramified Hecke character of infinity type $z^2$. Then there
exists a pair of characters $\phi=(\phi_1,\phi_2)$ satisfying
$(\phi)$ with $\chi=\phi_1/\phi_2$ such that  there is an
$\Ox$-algebra surjection
$$
\mathbf{T}/\mathbf{I}_{\phi} \twoheadrightarrow \Ox/ \left(
L^{\mathrm{int}}(0,\chi)  \right).
$$ \hspace{\fill}
\qedsymbol
\end{thm}

\subsection{Discussion of results} \label{discussion}
These are the first such cohomological congruences between
Eisenstein series and cuspforms for ${\rm GL}_2$ over an imaginary
quadratic field, except for the results for degree 2 Eisenstein
classes associated to unramified characters in \cite{F}.
 There are two options: either these congruences first arise for
${\rm GL}_{2/F}$ or they show that congruences over $\Q$ can be
lifted, in accordance with Langlands functoriality. The congruences
constructed in \cite{F} turn out to be base changes of congruences
over $\Q$ (see \cite{F} Satz 3.3).

Recall from \cite{GL} Theorem 2 and \cite{Cr} p. 413 that a cuspform
over $F$ is a base change if and only if its Hecke eigenvalues at
complex conjugate places coincide. Observe that if $\overline \chi
\neq \chi^c$ then the Hecke eigenvalues of our Eisenstein cohomology
class $\Eis$ (and all twists by a character) are distinct at complex
conjugate places (see (E2) for the definition of the eigenvalues).
Therefore in this case our congruences are new, i.e., are not base
changed.

If $\overline \chi = \chi^c$ then the proof of Lemma
\ref{anticycchar} implies that there exists a twist of the
Eisenstein class such that its eigenvalues at conjugate places
coincide. However, we cannot determine if the congruences are base
changed, as for cohomology in degree $1$ the arguments of \cite{HLR}
do not apply. We plan to investigate this question further.

\section{The case of unramified characters} \label{s6}
In this section we will prove Proposition \ref{thmunram}, i.e., show
the existence of an integral lift of the constant term of the
Eisenstein cohomology class $\Eis$, as defined in \S
\ref{Eisreview}. Our strategy is to find an involution on the
boundary cohomology such that the restriction map surjects onto the
$-1$-eigenspace of this involution, i.e., such that (for each
connected component of $\overline S_{K_f}$)
$$H^1(\Gamma \backslash\oh, \Ox)
\overset{\rm res}{\twoheadrightarrow} H^1(\partial(\Gamma\backslash
\oh),\Ox)^- \subset H^1(\partial(\Gamma\backslash \oh),\Ox), $$
where the superscript `-' indicates the $-1$-eigenspace. We prove
the existence of such an involution for all maximal arithmetic
subgroups of ${\rm SL}_2(F)$, extending a result of Serre for ${\rm
SL}_2(\mathcal{O})$. Proposition \ref{thmunram} is then proven by
showing that ${\rm res}(\Eis)$ lies in this $-1$-eigenspace.

\subsection{Involutions and the image of the restriction map}
Let $\Gamma \subset G(\Q)$ be an arithmetic subgroup. Given an
involution $\iota$ on $X=\Gamma \backslash \oh$ or $\partial(\Gamma
\backslash \oh)$ we define an involution on $H^1(X, R)$ via the
pullback of $\iota$ on the level of singular cocycles. Assuming that
we have an orientation-reversing involution on $\Gamma \backslash
\oh$ such that
$$H^1(\Gamma \backslash\oh, \Ox)
\overset{\rm res}{\rightarrow} H^1(\partial(\Gamma \backslash
\oh),\Ox)^- \subset H^1(\partial(\Gamma \backslash \oh),\Ox)$$ we
show that the map is, in fact, surjective. The existence of such an
involution will be shown for maximal arithmetic subgroups in the
following sections. We first recall:
\begin{thm}[(Poincar\'e and Lefschetz duality)]
Let $R$ be a Dedekind domain in which $2$ and $3$ are invertible.
Let $\iota$ be an orientation-reversing involution on $\Gamma
\backslash \oh$. Denoting by a superscript $+$ (resp. $-$) the
$+1$-(resp. $-1$-) eigenspaces for the induced involutions on
cohomology groups, we have perfect pairings
$$H^r_c(\Gamma \backslash \oh, R)^{\pm} \times H^{3-r}(\Gamma \backslash \oh, R)^{\mp} \to R \,
\text{ for } 0\leq r \leq 3$$ and
$$H^r(\partial(\Gamma \backslash \oh),R)^{\pm} \times H^{2-r}(\partial(\Gamma \backslash \oh),R)^{\mp} \to R \,
\text{ for } 0\leq r \leq 2.$$
Furthermore, the maps in the exact sequence
$$H^1(\Gamma \backslash\oh,R) \xrightarrow{\rm res} H^1(\partial(\Gamma \backslash \oh),R) \xrightarrow{\partial}
H^2_c(\Gamma \backslash\oh,R)$$ are adjoint, i.e.,
$$\langle {\rm res}(x),y\rangle=\langle x,\partial(y)\rangle.$$

\end{thm}

\begin{proof}[References]
Serre states this in the proof of Lemma 11 in \cite{Se} for field
coefficients, \cite{AS} Lemma 1.4.3 proves the perfectness for
fields R and \cite{U95} Theorem 1.6 for Dedekind domains as above.
Other references for this Lefschetz or ``relative" Poincar\'e
duality for oriented manifolds with boundary are \cite{Ma99}
Chapter 21, \S 4 and \cite{G67} (28.18). The pairings are given by
the cup product and evaluation on the respective fundamental
classes. We use that $\oh$ is an oriented manifold with boundary
and that $\Gamma$ acts on it properly discontinuously and without
reversing orientation. The lemma in \cite{F} \S 1.1 shows that the
order of any finite subgroup of $G(\Q)$ is divisible only by $2$
or $3$. See also \cite{Be} Theorem 5.1 and Lemma 5.2.
\end{proof}

\begin{lem} \label{surj}
Suppose in addition to the conditions of the previous theorem that
$R$ is a complete discrete valuation ring with finite residue field
of characteristic $p>2$. Suppose that we have an involution $\iota$
as in the Theorem such that
$$H^1(\Gamma \backslash \oh, R) \overset{\rm res}{\rightarrow}
H^1(\partial(\Gamma \backslash \oh),R)^{\epsilon},$$ where
$\epsilon=+1$ or $-1$. Then, in fact, the restriction map is
surjective.
\end{lem}

\begin{proof}
Let $\mathfrak{m}$ denote the maximal ideal of $R$. Since the
cohomology modules are finitely generated (so the Mittag-Leffler
condition is satisfied for $\underleftarrow{\lim} \,H^1( \cdot,
R/\mathfrak{m}^r)$), it suffices to prove the surjectivity for
each $r \in \N$ of
$$H^1(\Gamma \backslash \oh ,R/\mathfrak{m}^r) \twoheadrightarrow
H^1(\partial(\Gamma \backslash
\oh),R/\mathfrak{m}^r)^{\epsilon}.$$ For these coefficient systems
we are dealing with finite groups and can count the number of
elements in the image and the eigenspace of the involution; they
turn out to be the same.  We observe that $H^1(\partial(\Gamma
\backslash \oh),R/\mathfrak{m}^r)=H^1(\partial(\Gamma \backslash
\oh),R/\mathfrak{m}^r)^+ \oplus H^1(\partial(\Gamma \backslash
\oh),R/\mathfrak{m}^r)^-$ and that, by the last lemma, $$\#
H^1(\partial(\Gamma \backslash \oh),R/\mathfrak{m}^r)^+= \#
H^1(\partial(\Gamma \backslash \oh),R/\mathfrak{m}^r)^-.$$
Similarly we deduce from the adjointness of ${\rm res}$ and
$\partial$ and the perfectness of the pairings that ${\rm im}({\rm
res})^{\bot}= {\rm im}({\rm res})$ and so
$$\# {\rm im}({\rm res}) =\frac{1}{2} \# H^1(\partial(\Gamma
\backslash \oh),R/\mathfrak{m}^r).$$
\end{proof}

\subsection{Involutions for maximal arithmetic subgroups of ${\rm SL}_2(F)$}
For $\eta \in G(\Q)$ let $B^{\eta}$ be the parabolic subgroup
defined by $B^{\eta}(\Q)=\eta^{-1} B(\Q) \eta$. Let $\Gamma \subset
G(\Q)$ be an arithmetic subgroup. The set $\{B^{\eta}: [\eta] \in
B(\Q) \backslash G(\Q)/\Gamma \}$ is a set of representatives for
the $\Gamma$-conjugacy classes of Borel subgroups. Let $U^{\eta}$ be
the unipotent radical of $B^{\eta}$. For $D \in \mathbf{P}^1(F)$ let
$\Gamma_D=\Gamma \cap U_D$, where $U_D$ is the unipotent subgroup of
${\rm SL}_2(F)$ fixing $D$. Note that if $D_{\eta} \in
\mathbf{P}^1(F)$ corresponds to $[\eta] \in B(\Q) \backslash G(\Q) $
under the isomorphism of $B(\Q) \backslash G(\Q) \cong
\mathbf{P}^1(F)$ given by right action on $[0:1] \in
\mathbf{P}^1(F)$ then we have that $U_{D_{\eta}}=U^{\eta}(\Q)$ and
$\Gamma_{D_{\eta}}=\Gamma \cap U^{\eta}(\Q)=:\Gamma_{U^{\eta}}$.

Let $U(\Gamma)$ be the direct sum $\bigoplus_{[D] \in
\mathbf{P}^1(F)/\Gamma} \Gamma_{D}.$ Up to canonical isomorphism
this is independent of the choice of representatives $[D] \in
\mathbf{P}^1(F)/\Gamma$. The inclusions $\Gamma_{D} \hookrightarrow
\Gamma$ define a homomorphism $$\alpha: U(\Gamma) \to \Gamma^{\rm
ab}.$$

For $\Gamma={\rm SL}_2(\mathcal{O})$ \cite{Se} shows  that there is
a well-defined action of complex conjugation on $U({\rm
SL}_2(\mathcal{O}))$ induced by the complex conjugation action on
the matrix entries of $G_{\infty}={\rm GL_2}(\C)$. Denoting by $U^+$
the set of elements of $U({\rm SL}_2(\mathcal{O}))$ invariant under
the involution and by $U'$ the set of elements $u+ \overline u$ for
$u \in U({\rm SL}_2(\mathcal{O}))$, Serre proves:

\begin{thm}[(Serre \cite{Se} Th\'eor\`eme 9)] \label{TSerre}
For imaginary quadratic fields $F$ other than  $\Q(\sqrt{-1})$ or
$\Q(\sqrt{-3})$ the kernel of the homomorphism $\alpha:U({\rm
SL}_2(\mathcal{O})) \to {\rm SL}_2(\mathcal{O})^{\rm ab}$ satisfies
the inclusions $$6 U' \subseteq {\rm ker}(\alpha) \subseteq U^+.$$
\end{thm}
In the following we generalize this theorem to all maximal
arithmetic subgroups. After we had discovered this generalization we
found out that it had already been stated in \cite{BN}, but for our
application we need more detail than is provided there.

For $\mathfrak{b}$ a fractional ideal let
$$H(\mathfrak{b}):=\{\begin{pmatrix}a&b
\\ c&d
\end{pmatrix}\in {\rm SL}_2(F)|a, d \in \mathcal{O}, b \in \mathfrak{b}, c \in
\mathfrak{b}^{-1}\}.$$ This is a maximal arithmetic subgroup of
${\rm SL}_2(F)$ and any maximal arithmetic subgroup is conjugate to
$H(\mathfrak{b})$ (see \cite{EGM} Prop. 7.4.5). In order to study
the structure of $U(H(\mathfrak{b}))$ we define $j:\mathbf{P}^1(F)
\to \mathrm{Cl}(F)$ to be the map
$$j([z_1:z_2])=[z_1 \mathfrak{b} + z_2 \mathcal{O}].$$

\begin{thm} \label{bianchi}
For $\Gamma=H(\mathfrak{b})$, the induced map
$$j:\mathbf{P}^1(F)/\Gamma \to
\mathrm{Cl}(F)$$ is a bijection.
\end{thm}

\begin{proof}
Let $(x_1, x_2)$, $(y_1, y_2) \in F \times F$. It is easy to check
(see \cite{EGM} Theorem VII 2.4 for $\mathrm{SL}_2(\mathcal{O})$,
\cite{Be} Lemma 5.10 for the general case) that the following are
equivalent:

\emph{(1)} $x_1 \mathfrak{b} + x_2 \mathcal{O}  = y_1
\mathfrak{b}+y_2 \mathcal{O}.$

\emph{(2)} There exists $\sigma \in H(\mathfrak{b})$ such that
  $(x_{1}, x_{2})= (y_1,  y_2)\sigma.$

It remains to show the surjectivity of $j$. Given a class in ${\rm
Cl}(F)$ take $\mathfrak{a} \subset \mathcal{O}$ representing it.
By the Chinese Remainder Theorem one can choose $z_2 \in
\mathcal{O}$ such that
\begin{itemize}
  \item ${\rm ord}_{\wp}(z_2)={\rm
ord}_{\wp}(\mathfrak{a})$ if $\wp|\mathfrak{a}$.

  \item ${\rm ord}_{\wp}(z_2)=0$ if $\wp \nmid \mathfrak{a}, {\rm
  ord}_{\wp}(\mathfrak{b}) \neq 0.$
\end{itemize}
Then one chooses $z_1$ such that
\begin{itemize}
  \item ${\rm ord}_{\wp}(z_1 \mathfrak{b})>{\rm
ord}_{\wp}(z_2)$ if $\wp|\mathfrak{a}$ or ${\rm
ord}_{\wp}(\mathfrak{b}) \neq 0$.

  \item ${\rm ord}_{\wp}(z_1 \mathfrak{b})=0$ if $\wp |z_2, \wp
  \nmid \mathfrak{a}, $ and ${\rm ord}_{\wp}(\mathfrak{b})=0$.
\end{itemize}
These choices ensure that ${\rm ord}_{\wp}(z_1 \mathfrak{b} + z_2
\mathcal{O})={\rm ord}_{\wp}(\mathfrak{a})$ for all prime ideals
$\wp$.

\end{proof}

Following Serre \cite{Se} we now calculate explicitely
$\Gamma_{[z_1:z_2]}$ for $\Gamma=H(\mathfrak{b})$ and $[z_1:z_2]
\in \mathbf{P}^1(F)$.

\begin{lem}
For $\Gamma=H(\mathfrak{b})$, $\Gamma_{[z_1:z_2]}$ is conjugate in
$H(\mathfrak{b})$ to $$\{ \theta
\begin{pmatrix} 1&t \\ 0 & 1\end{pmatrix} \theta^{-1}  : t \in \mathfrak{a}^{-2} \mathfrak{b} \}
,$$ where $\mathfrak{a}=z_1 \mathfrak{b}+ z_2 \mathcal{O}$ and
$\theta$ is an isomorphism $\mathcal{O} \oplus \mathfrak{b}
\overset{\sim}{\to} \mathfrak{a} \oplus \mathfrak{a}^{-1}
\mathfrak{b}$ of determinant 1, i.e., such that its second
exterior power
$$ \Lambda^2 \theta: \Lambda^2 (\mathcal{O} \oplus \mathfrak{b})=
\mathfrak{b} \to \Lambda^2(\mathfrak{a} \oplus \mathfrak{a}^{-1}
\mathfrak{b}) = \mathfrak{a} \otimes \mathfrak{a}^{-1}
\mathfrak{b} =\mathfrak{b}$$ is the identity.

\end{lem}

\begin{proof}
The main change to \cite{Se} \S 3.6 is that we consider the
lattice $L:= \mathcal{O} \oplus \mathfrak{b}$ instead of
$\mathcal{O}^2$. We claim there exists a projective rank 1
submodule $E$ of $L$ containing a multiple of $(z_1, z_2)$. Let
$E$ be the kernel of the $\mathcal{O}$-homomorphism $
L=\mathcal{O} \oplus \mathfrak{b} \to F$ given by $(x,y) \mapsto y
z_1 -x z_2$. Since the image is $\mathfrak{a}=z_1 \mathfrak{b} +
z_2 \mathcal{O}$, we get $L/E \cong \mathfrak{a}$, so $L/E$ is
projective of rank 1 and $L$ decomposes as $E \oplus L/E$.

By definition $\Gamma_{[z_1:z_2]}$ fixes $L \cap \{\lambda (z_1,
z_2), \lambda \in F\}$, but this is exactly $E$. Since
$\Gamma_{[z_1:z_2]}$ is unipotent it can therefore be identified
with $\mathrm{Hom}_{\mathcal{O}}(L/E, E)$. For any fractional
ideal $\mathfrak{a}$, $\Lambda^2(\mathfrak{a})=0$ and so
$\mathfrak{b}=\Lambda^2(L)=\Lambda^2(E \oplus L/E)= E
\otimes_{\mathcal{O}} L/E$ so $E$ is isomorphic to $(L/E)^{-1}
\otimes \mathfrak{b}$. This implies an isomorphism
$\mathrm{Hom}_{\mathcal{O}}(L/E, E)= (L/E)^{-1} \otimes E \cong
(L/E)^{-1} \otimes (L/E)^{-1} \otimes \mathfrak{b} \cong
\mathfrak{a}^{-2} \mathfrak{b}$. Choosing an isomorphism $\theta:L
\to L/E \oplus E \cong \mathfrak{a} \oplus \mathfrak{a}^{-1}
\mathfrak{b}$ of determinant 1 we can represent
$\Gamma_{[z_1:z_2]}$ as stated above.
\end{proof}

Note that since $H(\mathfrak{b})$ is the stabilizer of any lattice
$\mathfrak{m} \oplus \mathfrak{n}$ with $\mathfrak{m}$ and
$\mathfrak{n}$ fractional ideals of $F$ such that
$\mathfrak{m}^{-1} \mathfrak{n} =\mathfrak{b}$, one can deduce

\begin{lem} \label{l5.7}
Let $\mathfrak{a}, \mathfrak{b}$ be two fractional ideals of $F$.
If $[\mathfrak{a}]=[\mathfrak{b}]$ in
$\mathrm{Cl}(F)/\mathrm{Cl}(F)^2$, then $H(\mathfrak{a}) =
H(\mathfrak{b})^{\gamma}$ with $\gamma \in \mathrm{GL}_2(F)$. If
the fractional ideals differ by the square of an
$\mathcal{O}$-ideal, then $\gamma$ can be taken to be in
$\mathrm{SL}_2(F)$.
\end{lem}

If the class of $\mathfrak{b}$ in $\mathrm{Cl}(F)$ is a square,
$H(\mathfrak{b})$ is isomorphic to $\mathrm{SL}_2(\mathcal{O})$ by
Lemma \ref{l5.7}, and the involution on
$U(\mathrm{SL}_2(\mathcal{O}))$ induced by complex conjugation and
Serre's Th\'eor\`eme 9 can easily be transferred to
$U(H(\mathfrak{b}))$. We therefore turn our attention to the case
when $$[\mathfrak{b}]  \text{ is not a square in } \mathrm{Cl}(F).$$ Note
that this implies that $[\mathfrak{b}]$ has even order, since any
odd order class can be written as a square.

Define an involution on $H(\mathfrak{b})$ to be the composition
of complex conjugation with an Atkin-Lehner involution, i.e., by

$$H=\begin{pmatrix} a&b\\c&d \end{pmatrix} \mapsto A \overline H
A^{-1}= \begin{pmatrix} \overline d& -{\rm Nm}(\mathfrak{b})
\overline c
\\- \overline b {\rm Nm}(\mathfrak{b})^{-1}& \overline{a}
\end{pmatrix},$$

where $A=\begin{pmatrix} 0&1\\-{\rm Nm}(\mathfrak{b})^{-1}&0
\end{pmatrix}$.

Like Serre, we will choose a set of representatives for the cusps
$\mathbf{P}^1(F)/H(\mathfrak{b})$ on which this involution acts.
For this we observe that if $\Gamma_{[z_1:z_2]}$ fixes $[z_1:z_2]$
then $A \overline \Gamma_{[z_1:z_2]} A^{-1}$ fixes $[\overline
z_1:\overline z_2] A^{-1}=[\overline z_2: -{\rm Nm}(\mathfrak{b})
\overline z_1]$. We use the isomorphism $j:
\mathbf{P}^1(F)/H(\mathfrak{b}) \to {\rm Cl}(F)$ to show that this
action on the cusps is fixpoint-free. We observe that if
$j([z_1:z_2])=\mathfrak{a}$ then $j([\overline z_1:\overline z_2]
A^{-1})=[\overline z_2 \mathfrak{b} + {\rm Nm}(\mathfrak{b})
\overline z_1 \mathcal{O}]=[\overline{\mathfrak{a}}
\mathfrak{b}]$. Note that $[\mathfrak{a}] \neq
[\overline{\mathfrak{a}} \mathfrak{b}]$ in ${\rm Cl}(F)$ since
otherwise $[\mathfrak{a}^2]=[{\rm Nm}(\mathfrak{a})
\mathfrak{b}]=[\mathfrak{b}]$, i.e., $[\mathfrak{b}]$ a square,
contradicting our hypothesis. So $\mathrm{Cl}(F)$ can be
partitioned into pairs $(\mathfrak{a}_i,\overline{\mathfrak{a}_i}
\mathfrak{b})$.

Choosing $[z_1^i:z_2^i]\in \mathbf{P}^1(F)$ such that
$\mathfrak{a}_i=z_1^i \mathfrak{b} + z_2^i \mathcal{O}$ we obtain
$$U(H(\mathfrak{b}))= \bigoplus_{(\mathfrak{a}_i,\overline{\mathfrak{a}_i}
\mathfrak{b})} (\Gamma_{[z_1^i:z_2^i]} \oplus A \overline
\Gamma_{[z_1^i:z_2^i]}  A^{-1} ).$$ Our choice of representatives
of $\mathbf{P}^1(F)/H(\mathfrak{b})$ shows that the involution
operates on $U(H(\mathfrak{b}))$ and, in fact, by identifying
$\Gamma_{[z_1^i:z_2^i]}$ with $ \{\theta
\begin{pmatrix} 1&s \\ 0 & 1 \end{pmatrix} \theta^{-1}: s \in
\mathfrak{a}_i^{-2} \mathfrak{b}  \}$ for $\theta: \mathcal{O}
\oplus \mathfrak{b}  \to \mathfrak{a}_i \oplus
\mathfrak{a}_i^{-1}\mathfrak{b}$ and $A \overline
\Gamma_{[z_1^i:z_2^i]} A^{-1}$
with $\{\theta'  \begin{pmatrix} 1&0 \\
-t & 1 \end{pmatrix} \theta'^{-1}:t \in
\overline{\mathfrak{a}_i}^{-2} \mathfrak{b}^{-1} \}$ for
$\theta'=A \overline{\theta} A^{-1} : \mathcal{O} \oplus
\mathfrak{b} \to \overline{\mathfrak{a}_i}^{-1} \oplus
\overline{\mathfrak{a}_i}\mathfrak{b} $, we can describe the
involution on each of the pairs as

$$(s, t) \in \mathfrak{a}_i^{-2} \mathfrak{b} \oplus
\overline{\mathfrak{a}_i}^{-2} \mathfrak{b}^{-1} \mapsto
(\overline{t} {\rm Nm}(\mathfrak{b}), \overline s {\rm
Nm}(\mathfrak{b})^{-1}).$$

Now denote by $U^+$ the set of elements of $U(H(\mathfrak{b}))$
invariant under the involution $H \mapsto A \overline H A^{-1}$,
and by $U'$ the set of elements $u+A \overline{u}A^{-1}$ for $u
\in U(H(\mathfrak{b}))$.

\begin{thm} \label{genserre}
For $\Gamma=H(\mathfrak{b})$ with $[\mathfrak{b}]$ a non-square in
$\mathrm{Cl}(F)$, the kernel $N$ of the homomorphism
$$\alpha :U(\Gamma)\to \Gamma^{\mathrm{ab}}$$
coming from the inclusion $\Gamma_D \hookrightarrow \Gamma$ for $D
\in \mathbf{P}^1(F)$ satisfies $6U' \subset N \subset U^+$.
\end{thm}

\begin{proof}
With small modifications, we follow Serre's proof of his
Th\'{e}or\`{e}me 9. As in Serre's case, it suffices to prove the
inclusion $6U' \subset N$, i.e., that $6(u+ A \overline{u}A^{-1})$
maps to an element of the commutator
$[H(\mathfrak{b}),H(\mathfrak{b})]$:

Suppose that we have $6U' \subset N$, but that there exists an
element $u \in N$ not contained in $U^+$. Then the subgroup of $N$
generated by $6U'$ and $u$ has rank $\#{\rm Cl}(F) +1$. This
contradicts the fact that the kernel of $\alpha$ has rank $\#{\rm
Cl}(F)$ (see \cite{Se} Th\'eor\`eme 7). (The latter is proven by
showing dually that the rank of the image of the restriction map
$H^1(H(\mathfrak{b}) \backslash \oh, R) \to
H^1(\partial(H(\mathfrak{b}) \backslash \oh ), R)$ has half the
rank of that of the boundary cohomology. This we showed in the
proof of Lemma \ref{surj}).

To prove $6 U' \subset N$ we make use of Serre's Proposition 6:

\begin{prop}[(\cite{Se} Proposition 6)]
Let $\mathfrak{q}$ be a fractional ideal of $F$ and let $t\in
\mathfrak{q}$ and $t'=\overline{t}/{\rm Nm}(\mathfrak{q})$ so that
$t' \in \mathfrak{q}^{-1}$. Put $x_t=\begin{pmatrix} 1&t
\\ 0 & 1 \end{pmatrix}$ and $y_t=\begin{pmatrix} 1&0 \\ -t' & 1
\end{pmatrix}$. Then $(x_t y_t)^6$ lies in the commutator subgroup
of $H(\mathfrak{q})$.
\end{prop}

Put $\mathfrak{a}:= z_1 \mathfrak{b} + z_2 \mathcal{O}$.
If $u \in \Gamma_{[z_1:z_2]}$,  identify it with $\theta^{-1}\begin{pmatrix} 1&t \\
0 & 1
\end{pmatrix} \theta$ for some $t \in \mathfrak{a}^{-2}
\mathfrak{b}$ and $\theta: \mathcal{O} \oplus \mathfrak{b} \to
\mathfrak{a} \oplus \mathfrak{a}^{-1}\mathfrak{b}$ of determinant
1. One easily checks that $A \overline u A^{-1}$ then corresponds
to $(A \overline{\theta} A^{-1}) \begin{pmatrix} 1&0 \\
-\overline{t} {\rm Nm}(\mathfrak{b})^{-1} & 1
\end{pmatrix} (A \overline{\theta} A^{-1})$. Like Serre, we use
that since $[\overline{\mathfrak{a}}]=[\mathfrak{a}^{-1}]$, $A
\overline u A^{-1}$ is also given by Theorem \ref{bianchi}
by $B^{-1} \theta^{-1} \begin{pmatrix} 1&0 \\
-t'  & 1 \end{pmatrix} \theta B$ for $t'=\overline{t} {\rm
Nm}(\mathfrak{b})^{-1} {\rm Nm}(\mathfrak{a})^2$ and $B \in
H(\mathfrak{b})$ taking $\begin{pmatrix} {\rm Nm}(\mathfrak{b})
\overline z_2\\ \overline z_1 \end{pmatrix}$ to ${\rm
Nm}(\mathfrak{a})^{-1}
\begin{pmatrix} {\rm Nm}(\mathfrak{b}) \overline z_2\\ \overline z_1
\end{pmatrix}$.

  Since
$\theta^{-1} x_t y_t \theta$ is a representative of $u + B A
\overline u A^{-1} B^{-1}$, we deduce from the above Proposition
with $\mathfrak{q}=\mathfrak{a}^{-2} \mathfrak{b}$ that $6(u + B A
\overline u A^{-1} B^{-1})$ and therefore $6(u+ A \overline u
A^{-1})$ lie in $[H(\mathfrak{b}), H(\mathfrak{b})]$.
\end{proof}

The following observation links $U(\Gamma)$ to the cohomology of the
boundary components:

\begin{lem} \label{l5.6}
For imaginary quadratic fields $F$ other than  $\Q(\sqrt{-1})$ or
$\Q(\sqrt{-3})$, $\Gamma \subset {\rm SL}_2(F)$ an arithmetic
subgroup, $P$ a parabolic subgroup of ${\rm Res}_{F/\Q}({\rm
SL}_{2/F})$ with unipotent radical $U_P$, and $R$ a ring in which 2
is invertible we have
$$H^1(\Gamma_{P}, R) \cong H^1(\Gamma_{U_P}, R), $$ where
$\Gamma_P=\Gamma \cap P(\Q)$ and $\Gamma_{U_P}=\Gamma \cap U_P(\Q)$.
\end{lem}

\begin{proof}
Serre shows in \cite{Se} Lemme 7 that $\Gamma_{U_P} \triangleleft
\Gamma_{P}$ and that the quotient $W_P=\Gamma_P/\Gamma_{U_P}$ can be
identified with a subgroup of the roots of unity of $F$, i.e., of
$\{ \pm 1\}$ since $F\neq \Q(\sqrt{-1}), \Q(\sqrt{-3})$. The Lemma
follows from the Inflation-Restriction sequence. See also \cite{Ta}
p.110.
\end{proof}

By (\ref{bdry}), (\ref{shgp}), and Lemma \ref{l5.6} we have
\begin{equation} \label{bdgp} H^1(\partial(\Gamma \backslash \oh),R) \cong \coprod_{[\eta]
\in \mathbf{P}^1(F)/\Gamma} H^1(\Gamma_{U^{\eta}},R)= H^1(U(\Gamma),
R).\end{equation}
 We now reinterpret Serre's Theorem and its generalization as follows:

\begin{prop} \label{involutions}
For imaginary quadratic fields $F$ other than  $\Q(\sqrt{-1})$ or
$\Q(\sqrt{-3})$ and $R$ a ring in which 2 and 3 is invertible, the
image of the restriction map
$$H^1(\Gamma \backslash \oh,R) \overset{{\rm res}}{\to} H^1(\partial(\Gamma \backslash \oh),R)$$
is contained in the $-1$-eigenspace of the involution induced by

\begin{itemize}
\item    $\iota:\h \to \h: (z, t)  \mapsto
(\overline z, t)$  if  $\Gamma={\rm SL}_2(\mathcal{O})$
\item $\iota: \h \to \h: (z, t) \mapsto A.(\overline z, t)$
for $A=\begin{pmatrix} 0&1\\-{\rm Nm}(\mathfrak{b})^{-1}&0
\end{pmatrix}$ if $\Gamma=H(\mathfrak{b})$ with
$[\mathfrak{b}]$ a non-square in $\mathrm{Cl}(F)$.
  \end{itemize} and these involutions are orientation-reversing.
\end{prop}

By Lemma \ref{surj} this immediately implies:
\begin{cor} \label{torsion2}
For imaginary quadratic fields $F$ other than  $\Q(\sqrt{-1})$ or
$\Q(\sqrt{-3})$, $\Gamma={\rm SL}_2(\mathcal{O})$ or $
H(\mathfrak{b})$ with $[\mathfrak{b}]$ a non-square in
$\mathrm{Cl}(F)$, and $R$ a complete discrete valuation ring in
which 2 and 3 are invertible and with finite residue field of
characteristic $p>2$, the restriction map
$$H^1(\Gamma \backslash \oh,R) \overset{\rm res}{\to} H^1(\partial(\Gamma \backslash \oh),R)^{-}$$
surjects onto the $-1$-eigenspace of the involutions defined in the
Proposition.
\end{cor}

\begin{proof}[Proof of Proposition]
Write $I:\Gamma \to \Gamma$ for the involution $$
  \begin{cases}
    \gamma \mapsto \overline \gamma & \text{if } \Gamma={\rm SL}_2(\mathcal{O}), \\
    \gamma \mapsto A \overline \gamma A^{-1} & \text{if } \Gamma=H(\mathfrak{b}).
  \end{cases}$$

The involutions $\iota$ extend canonically to $\oh$. One checks
that for $\gamma \in \Gamma$ we have
\begin{equation} \label{comp} \iota(\gamma.(z,t))= I(\gamma)\iota(z,t). \end{equation} This implies
that the involutions operate on $\Gamma \backslash \h$ and $\Gamma
\backslash \oh$, and hence on $\partial(\Gamma \backslash \oh)$.
To show that they act by reversing the orientation note that
complex conjugation corresponds to reflection in a half-plane of
$\h$ and therefore reverses the orientation.  Furthermore, ${\rm
GL}_2(\C)$ acts on $\h$ via $A'=({\rm det}(A)^{-\frac{1}{2}}) A
\in {\rm SL}_2(\C)$ and ${\rm SL}_2(\C)$ acts without reversing
orientation, as can be seen from the geometric definition of its
action via the Poincar\'e extension of the action on
$\mathbf{P}^1(\C)$ (see \cite{EGM} pp.2-3).

Using (\ref{comp}) one shows that under the isomorphism
$$H^1(\partial(\Gamma \backslash \oh),R) \overset{(\ref{bdgp})}{\cong} H^1(U(\Gamma), R)$$
$\iota$ corresponds to the involution on $H^1(U(\Gamma), R)={\rm
Hom}(U(\Gamma), R)$ given by $\varphi \mapsto I(\varphi)$, where
$I(\varphi)(u):=\varphi(I(u))$.

We can therefore check that the image of the restriction maps is
contained in the $-1$-eigenspace  on the level of group
cohomology: The restriction map is given by
$${\rm Hom}(\Gamma^{\rm ab}, R) \to {\rm Hom}(U(\Gamma), R):
\varphi \mapsto \varphi \circ \alpha.$$ By Serre's theorem and
Theorem \ref{genserre}, $0=\varphi(\alpha(u I(u)))
=\varphi(\alpha(u)) + \varphi(\alpha(I(u)))$, so $I(\varphi \circ
\alpha)(u)=\varphi(\alpha(I(u)))=-\varphi(\alpha(u))$ for any $u
\in U(\Gamma)$.
\end{proof}

\subsection{Integral lift of constant term}
Recall the statement of Proposition \ref{thmunram}:
  \begin{prop}[(=Proposition \ref{thmunram})]
Let $\chi=\phi_1/\phi_2$ be an unramified Hecke character of
infinity type $z^2$. Assume that $1$ is the only unit in
$\mathcal{O}^*$ congruent to $1$ modulo the conductor of $\phi_1$.
Let $K_f$ and $\Eis$ be defined as in Section \ref{Eisreview}.
Assume  ${\rm
  res}(\Eis) \in \widetilde H^1(\partial \overline S_{K_f}, \Ox)$.
  Then there exists $c_{\phi}
  \in \widetilde H^1(S_{K_f},\Ox)$ with $${\rm res}(c_{\phi})={\rm
  res}(\Eis) \in \widetilde H^1(\partial \overline S_{K_f}, \Ox).$$
  \end{prop}

First observe that everywhere unramified characters with infinity
type $z^2$ exist only for $F \neq \Q(\sqrt{-1}), \Q(\sqrt{-3})$.
For unramified $\chi$ we have
$$K_f= \prod_{v \mid \mathfrak{M}_1} U^1(\mathfrak{M}_{1,v}) \prod_{v \nmid \mathfrak{M}_1}
{\rm GL}_2(\mathcal{O}_v).$$  Recall that
$U^1(\mathfrak{M}_{1,v})=\{k \in {\rm GL}_2(\mathcal{O}_v) : {\rm
det}(k) \equiv 1 \mod{\mathfrak{M}_{1,v}} \}$. By assumption, $1$
is the only unit in $\mathcal{O}^*$ congruent to $1$ modulo
$\mathfrak{M}_1$ so we get $K_f \cap {\rm GL}_2(F)={\rm
SL}_2(\mathcal{O})$.

Recall from \S \ref{s2.5} the decomposition of $S_{K_f}$ into its
connected components. The above implies that we can write
$S_{K_f}$ as a disjoint union of $\Gamma \backslash \h$ with
$\Gamma=H(\mathfrak{b})$ for suitable fractional ideals
$\mathfrak{b}$: For a finite idele $a$, denote by $(a)$ the
corresponding fractional ideal. Write
$$S_{K_f} \cong \coprod_{i=1}^{\# \pi_0(K_f)} \Gamma_{t_i}\backslash
\h ,$$ where $\Gamma_{t_i}=G(\Q) \cap t_i K_f t_i^{-1}$ and the
$t_i \in G(\A_f)$ are given by $t_i=\begin{pmatrix} \gamma_j a_k
b_m&0\\0&b_m
\end{pmatrix}$, with
\begin{itemize}
  \item $\{ \gamma_j \}$ a system of representatives  of
$$\mathrm{ker}(\pi_0(K_f) \to \mathrm{Cl}(F)) \cong \mathcal{O}^*  \backslash \prod_v
\mathcal{O}_v^* / \mathrm{det}(K_f),$$
  \item $\{a_k \}$ a set of representatives  of
 $\mathrm{Cl}(F)/(\mathrm{Cl}(F))^2$ in $\A_{F,f}^*$ (and we represent the principal ideals by
 $(1)$),
  \item $\{b_m^2 \}$ a set representing $\mathrm{Cl}(F)^2$.
\end{itemize}

Note that for these choices $\Gamma_{t_i}=H((a_k))$ and either
$a_k=1$ or $[(a_k)]$ is not a square in $\mathrm{Cl}(F)$.
 This allows us to apply our
results for maximal arithmetic subgroups from the previous section
by considering the restriction maps to the boundary separately for
each connected component.

\begin{prop}
$$[{\rm res}(\Eis)] \in
(H^1(\partial \overline S_{K_f}, \Ox)^-)_{\rm free}, $$ where
$H^1(\partial \overline S_{K_f}, \Ox)^-$ is defined via the
isomorphism to
$$\bigoplus_{i=1}^{\#\pi_0(K_f)}  H^1(\partial(\Gamma_{t_i} \backslash \oh),
\Ox)^- $$ where `-' indicates the $-1$-eigenspace of the
involutions defined in Proposition \ref{involutions}.
\end{prop}

\begin{rem}
Together with Corollary \ref{torsion2} this shows the existence of
an integral lift of the constant term and proves Proposition
\ref{thmunram}.
\end{rem}

\begin{proof}
We will consider the restriction maps to the boundary separately
for each connected component $\Gamma_{t_i} \backslash \h$:
$$H^1(\Gamma_{t_i} \backslash\h, \Ox)
\overset{\rm res}{\to} H^1(\partial(\Gamma_{t_i} \backslash \oh),
\Ox ) \overset{(\ref{bdry})}{\cong} \bigoplus_{[\eta] \in
\mathbf{P}^1(F)/\Gamma_{t_i}} H^1(\Gamma_{t_i, B^{\eta}} \backslash
\h, \Ox),$$ where $\Gamma_{t_i, B^{\eta}}= \Gamma_{t_i} \cap
\eta^{-1} B(\Q) \eta.$ By (\ref{shgp}) and Lemma \ref{l5.6} we have
$ H^1(\Gamma_{t_i, B^{\eta}} \backslash \h, \Ox) \cong
H^1(\Gamma_{t_i, U^{\eta}}, \Ox)$. We recall from \cite{Be06}
Proposition 10, Lemma 11, and Proof of Proposition 16 that ${\rm
res}(\Eis)$ restricted to this boundary component is represented by
\begin{equation} \label{gpconst} \eta_{\infty}^{-1} \begin{pmatrix} 1&x \\0&1 \end{pmatrix}
\eta_{\infty} \mapsto x \Psi_{\phi}(\eta_f t_i) -
\frac{L(0,\overline \chi)}{L(0,\chi)} W(\chi) \cdot \overline x
\Psi_{w_0.\phi}(\eta_f t_i),\end{equation} where $W(\chi)$ is the
root number of $\chi$, $\eta_f$ and $\eta_{\infty}$ denote the
images of $\eta \in G(\Q)$ in $G(\A_f)$ and $G(\Real)$,
respectively, $w_0.(\phi_1,\phi_2)=(\phi_2 \cdot | \cdot |, \phi_1
\cdot | \cdot |^{-1})$, and $\Psi_{\phi}:G(\A_f) \to \C$ satisfies
$$\Psi_{\phi}(\begin{pmatrix} a&b\\0&d
\end{pmatrix}  k)=\phi_1(a) \phi_2(d) \text{ for }
\begin{pmatrix}a&b\\0&d
\end{pmatrix}\in B(\A_f), k \in \prod_v {\rm
SL}_2(\mathcal{O}_v)\subset K_f.$$ Note that, in particular,
$\Psi_{\phi}(bg)=\phi_{\infty}^{-1}(b)\Psi_{\phi}(g)$ for $b\in
B(F)\subset G(\A_f)$, where we use the convention introduced in \S
\ref{s2.5} for considering $\phi$ as a character of $B(\A)$. By
Lemma \ref{unanticyc}, $\chi^c=\overline \chi$, so
$L(0,\chi)=L(0,\overline \chi)$. Furthermore, $W(\chi)=i^2 (\chi/|
\chi |)(\delta_F^{-1})$ for $\delta_F$ the different ideal of
$F/\Q$ (see e.g. the proof of Proposition 2.1.7 of \cite{AH}).
Since $\delta_F^{-1}=g'(\alpha) \mathcal{O}$, where
$\mathcal{O}=\Z[\alpha]$ and $g \in \Z[X]$ is the minimal
polynomial of $\alpha$ over $\Q$ (see \cite{FT} III (2.20)) one
checks that for all imaginary quadratic fields there exists a
generator $\delta$ of $\delta_F$ satisfying $\overline
\delta=-\delta$. We deduce therefore that in our case $W(\chi)=1$.

We need to prove that (\ref{gpconst}) lies in the $-1$-eigenspace
of the involution induced by $u \mapsto \overline u$ for
$\Gamma_{t_i}={\rm SL}_2(\mathcal{O})$ and by $u \mapsto A
\overline u A^{-1}$ for $\Gamma_{t_i}=H(\mathfrak{b})$, where
$A=\begin{pmatrix} 0&1\\-N^{-1}&0 \end{pmatrix}$ with $N={\rm
Nm}(\mathfrak{b})$.

Case $\Gamma_{t_i}={\rm SL}_2(\mathcal{O})$: Recall that in this
case $t_i=\begin{pmatrix}\gamma_i b_i&0\\0&b_i
\end{pmatrix}$ for some $\gamma_i \in \hat{\mathcal{O}}^*$ and $b_i \in \A_{F,f}^*$.

It suffices to prove that $\Psi_{\phi}(\eta_f t_i
)=\Psi_{w_0.\phi}(\overline \eta_f t_i)$. For this we use the
Bruhat decomposition of matrices in ${\rm GL}_2(F)$ given by:
$$
\begin{pmatrix} a&b\\c&d
\end{pmatrix}=
  \begin{cases}
    \begin{pmatrix}1&b/d\\0&1 \end{pmatrix} \begin{pmatrix} a&0\\0&d\end{pmatrix}& \text{ if } c=0, \\
    \begin{pmatrix} 1&a/c\\0&1\end{pmatrix} \begin{pmatrix} \frac{ad-bc}{c}&0\\0&-c\end{pmatrix}
    \begin{pmatrix} 0&1\\-1&0\end{pmatrix}
    \begin{pmatrix}1&d/c\\0&1 \end{pmatrix}& \text{ otherwise}.
  \end{cases}
$$

Since $\Psi_{\phi}(\begin{pmatrix} a&b\\0&d
\end{pmatrix}g)=\Psi_{\phi}(\begin{pmatrix} a&b\\0&d
\end{pmatrix}) \Psi_{\phi}(g)$ we can consider separately the cases
\begin{enumerate}[(a)]
  \item $\eta=\begin{pmatrix} a&b\\0&d
\end{pmatrix}$ for $a, b, d \in F$ and
  \item $\eta=\begin{pmatrix} 0&1\\-1&0\end{pmatrix}\begin{pmatrix}
  1&e\\0&1\end{pmatrix}$ for $e\in F$.
\end{enumerate}

We check that for (a) $ \Psi_{\phi}(\eta_f
\begin{pmatrix}\gamma_i b_i&0\\0&b_i
\end{pmatrix})= \phi_1(\gamma_i b_i)
\phi_2(b_i)\Psi_{\phi}(\eta_f)$ and $\Psi_{w_0.\phi}(\overline
\eta_f \begin{pmatrix}\gamma_i b_i&0\\0&b_i
\end{pmatrix})=\phi_2(\gamma_i)|\gamma_i| \phi_1(b_i) \phi_2(b_i) \Psi_{w_0.\phi}(\overline \eta_f). $
Since $\gamma_i \in \hat{\mathcal{O}}^*$ and $\chi=\phi_1/\phi_2$
is unramified it suffices to show that $\Psi_{\phi}(\eta_f
)=\Psi_{w_0.\phi}(\overline \eta_f)$. In case (b) we similarly
reduce to this assertion.

In (a) we get
$\Psi_{\phi}(\eta_f)=\phi_{1,\infty}^{-1}(a)\phi_{2,\infty}^{-1}(d)=\frac{d}{a}$.
Since $w_0.\phi$ has infinity type $(\overline z, \overline
z^{-1})$ this equals $\Psi_{w_0.\phi}(\overline \eta_f)$. In (b)
we need to calculate the Iwasawa decomposition of $\eta$ in ${\rm
GL}_2(F_v)$ if $e \notin \mathcal{O}_v$ (at all other places
$\Psi_{\phi}(\eta_v)=\Psi_{w_0.\phi}(\overline \eta_{\overline v}
)=1$). It is given by
$$\begin{pmatrix} 0&1\\-1&0\end{pmatrix}\begin{pmatrix}
  1&e\\0&1\end{pmatrix}=\begin{pmatrix} e^{-1}&0\\0&e\end{pmatrix}\begin{pmatrix}
  -1&0\\-e^{-1}&-1\end{pmatrix}.$$
So, if $e \notin \mathcal{O}_v$ then
$\Psi_{\phi}(\eta_v)=(\phi_2/\phi_1)_v(e)=\chi^{-1}_v(e)$, which
matches $\Psi_{w_0.\phi}(\overline \eta_{\overline
v})=(\phi_1/\phi_2)_{\overline v}(\overline e) |\overline
e|_{\overline v}^{-2}$ because $\chi^c=\overline \chi$ and $\chi
\overline \chi=| \cdot |^2$.

Case $\Gamma_{t_i}=H(\mathfrak{b})$: The involution maps the cusp
corresponding to $B^{\eta}$ to $B^{\overline \eta A^{-1}}$. We
therefore have to prove that
\begin{equation}\label{Psi=}\Psi_{\phi}(\eta_f t_i)=\Psi_{w_0.\phi}(\overline \eta_f
A^{-1} t_i).\end{equation} Recall that $t_i=\begin{pmatrix}x_i
b_i&0\\0&b_i \end{pmatrix}$ for some $x_i, b_i \in \A_{F,f}^*$.
Again making use of the Bruhat decomposition, we need to only
consider $\eta$ as in cases (a) and (b) above. Following the
arguments used for Case (1), Case(a) reduces immediately to
showing that $\Psi_{\phi}(t_i)=\Psi_{w_0.\phi}(A^{-1} t_i).$ The
left hand side equals $\phi_{1,f}(x_i b_i) \phi_{2,f}(b_i)$, the
right hand side is
\begin{eqnarray*} \Psi_{w_0.\phi}(\begin{pmatrix} N&0\\0&1
\end{pmatrix}\begin{pmatrix}
0&1\\-1&0 \end{pmatrix}\begin{pmatrix} x_i b_i&0\\0&b_i
\end{pmatrix}) &=&N^{-1}\Psi_{w_0.\phi}(\begin{pmatrix} b_i&0\\0&x_i
b_i
\end{pmatrix})\\&=&N^{-1} \phi_{1,f}(x_i b_i) \phi_{2,f}(b_i)
|x_i|_f^{-1}.\end{eqnarray*} Equality follows from
$|x_i|_f^{-1}={\rm Nm}(\mathfrak{b})$.

For (b), one quickly checks that for $\eta=\begin{pmatrix}
0&1\\-1&0\end{pmatrix}$ the two sides in (\ref{Psi=}) agree. For
general $\eta=\begin{pmatrix}
0&1\\-1&0\end{pmatrix}\begin{pmatrix}
  1&e\\0&1\end{pmatrix}$ one shows that, on the one hand,
  $$\eta_f \begin{pmatrix} x_i b_i&0\\0&b_i
\end{pmatrix}= \begin{pmatrix} b_i&0\\0&x_i b_i \end{pmatrix}\begin{pmatrix}
0&1\\-1&0\end{pmatrix}\begin{pmatrix}
  1&e x_i\\0&1\end{pmatrix},$$
  and on the other hand,
  $$\overline \eta_f A^{-1} \begin{pmatrix} x_i b_i&0\\0&b_i
\end{pmatrix}=\begin{pmatrix} x_i b_i&0\\0&b_i N
\end{pmatrix}\begin{pmatrix}
0&1\\-1&0\end{pmatrix}\begin{pmatrix}
  1&\overline e x_i/N\\0&1\end{pmatrix}.$$

Since $(x_i \overline x_i)=(N)$ the valuations of $\overline e
x_i/N$ agrees with that of $\overline e \overline x_i$. Repeating
the calculation for $\eta=\begin{pmatrix} 0&1\\-1&0\end{pmatrix}$
and then applying the argument from Case 1(b) (since $\chi$ is
unramified we are only concerned about the valuation of the upper
right hand entry) we also obtain equality.

\end{proof}

\section{Bloch-Kato Conjecture for Hecke characters} \label{s7}
Combining Theorem \ref{Selthm} with Theorem \ref{condthm} or
\ref{unramcor} we get lower bounds on the size of the Selmer group
of Hecke characters. We want to demonstrate this application under
the particular conditions of Theorem \ref{unramcor} and relate it to
the Bloch-Kato conjecture.
\begin{thm} \label{smallthm}
Assume in addition that $p>3$, $p \nmid \#{\rm Cl}(F)$ and that
$\ell \not \equiv \pm 1 \mod{p}$ for $\ell \mid d_F$. If $p$ is
inert in $F/\Q$ then assume Conjecture \ref{chrys}. Let $\chi$ be an
unramified Hecke character of infinity type $z^2$. Then we have
$$\mathrm{val}_p \#{\rm Sel}(F,\chi_{\mathfrak{p}} \epsilon) \geq
\mathrm{val}_p(\# \Ox/L^{\mathrm{int}}(0,\chi))$$
\end{thm}

\begin{proof}
Put $\rho:=\chi_{\mathfrak{p}} \epsilon$. Theorem \ref{Selthm}
together with Lemma \ref{anticycchar} and Theorem \ref{unramcor}
imply
$$\mathrm{val}_p \#{\rm Sel}^{\Sigma_{\phi}\backslash \Sigma_{\rho}}(F,\chi_{\mathfrak{p}} \epsilon)
\geq \mathrm{val}_p(\# \Ox/L^{\mathrm{int}}(0,\chi)),$$ where
$\Sigma_{\rho}=\{v \mid p\}$ and $\phi=(\phi_1, \phi_2)$ is given
by Theorem \ref{unramcor}. For the definition of $\Sigma_{\phi}$
see the start of \S \ref{SelEis}. Recall that by $(\phi)$ the set
$\Sigma_{\phi}\backslash \{v \mid p\}$ contains only places $v$
such that $\overline v=v$ and $\# \mathcal{O}_v/\mathfrak{P}_v
\not \equiv \pm 1 \mod{p}$. By Lemma \ref{unanticyc} we have
$\chi^c=\overline \chi$, which implies that $\rho$ is
anticyclotomic, and so we get $\rho({\rm Frob}_v)=\rho({\rm
Frob}_v^c)=\rho^{-1}({\rm Frob}_v)$, or $\rho({\rm Frob}_v)= \pm
1.$ Hence we have ensured that
$$\rho({\rm Frob}_v) \not \equiv \epsilon({\rm Frob}_v) \mod{p}$$
 for all $v \in \Sigma_{\phi}\backslash \Sigma_{\rho}$, so the theorem follows
 from applying Lemma \ref{tameinertia}.
\end{proof}

\begin{example}
A numerical example in which the conditions of our Theorem are
satisfied and a non-trivial lower bound on a Selmer group is
obtained is given by the following: Let $F=\Q(\sqrt{-67})$ and
$p=19$. One checks that $19$ splits in $F$. Since the class number
is $1$, there is only one unramified Hecke character of infinity
type $z^2$. Up to $p$-adic units $L^{\mathrm{alg}}(0,\chi)$ is
given by $\frac{L(0,\chi)}{\Omega^2}$ where $\Omega$ is the Neron
period of the elliptic curve $y^2+y=x^3- 7370x+243582$, which has
conductor $67^2$ and complex multiplication by $\mathcal{O}$.
Using MAGMA and ComputeL \cite{ComputeL} one calculates that
$L^{\mathrm{alg}}(0,\chi) \in \Z_{19}$ and
$$\mathrm{val}_{19}(L^{\mathrm{alg}}(0,\chi))=1.$$
\end{example}

\subsection{Comparison with other results}
Assume from now on that $\#{\rm Cl}(F)=1$. Let $\Psi:F^* \backslash
\A_F^* \to \C^*$ be a Hecke character of infinity type $z^{-1}$
which satisfies $\Psi^c=\overline \Psi$. Then there exists an
elliptic curve $E$ over $\Q$ with complex multiplication by
$\mathcal{O}$ and associated Gr\"ossencharacter $\Psi$. Consider
$$\rho=(\Psi^{k} \overline{\Psi^{-j}})_{\mathfrak{p}} \text{ for }
k>0, j \geq 0.$$

We now have the following proposition from \cite{Dee} (Proposition
4.4.3 and \S5.3):
\begin{prop} [(Dee)]
The group ${\rm Sel}(F,\rho)$ is finite if and only if ${\rm
Sel}(F,\rho^{-1} \epsilon)$ is finite. If this is the case then
$$\#{\rm Sel}(F,\rho)=\#{\rm Sel}(F,\rho^{-1} \epsilon)< \infty.$$
\end{prop}
By considering $\chi=\Psi^{-2}$ for some $\Psi$ as above, compare
therefore Theorem \ref{smallthm} with the following:
\begin{thm}[(Han, \cite{Han})] Suppose $k>j+1$. For inert $p$
also assume that $\rho$ is non-trivial when restricted to ${\rm
Gal}(F(E_{\mathfrak{p}})/F)$. Then ${\rm Sel}(F,\rho)$ is finite
and
$$\mathrm{val}_p \# {\rm Sel}(F,\rho)=\mathrm{val}_p(\# \Ox/L^{\rm
int}(0, \Psi^{-k} \overline{\Psi}^{j})).$$
\end{thm}

Previously, Kato proved this in the case $k>0$ and $j=0$, cf.
\cite{Kato}. For a similar result in the case of split $p$ see
\cite{Guo93}. Han claims that his method extends to general class
numbers. All proofs take as input the proof of the Main Conjecture
of Iwasawa theory by Rubin \cite{Ru}.

We refer to \cite{Guo96} \S 3 for the proof that this statement on
the size of the Selmer group is equivalent to the (critical case of
the) $p$-part of the Bloch-Kato Tamagawa number conjecture for the
motives associated to the Hecke characters.

\begin{acknowledgements}
The results of this paper generalize part of my thesis \cite{Be}
with Chris Skinner at the University of Michigan. The author would
like to thank Trevor Arnold, John Cremona, Kris Klosin, Chris
Skinner, Jacques Tilouine, and Eric Urban for helpful discussions.
Thanks also to the Max Planck Institute in Bonn for their
hospitality and support during the writing of this article.
\end{acknowledgements}

\bibliographystyle{amsalpha}
\bibliography{biblio3}
\end{document}